\documentclass[12pt]{amsart}
\usepackage{amssymb,graphicx}
\usepackage{lmodern}
\usepackage{tikz}
\usetikzlibrary{arrows.meta,fit,shapes.geometric}
\usepackage{enumitem}
\usepackage[margin=2.5cm]{geometry}
\usepackage{mathrsfs}

\newtheorem{theorem}{Theorem}
\newtheorem{corollary}[theorem]{Corollary}
\newtheorem{lemma}[theorem]{Lemma}

\theoremstyle{definition}

\title{The Iterated Local Transitivity Model for Oriented Graphs}

\author[A.\ Bonato]{Anthony Bonato${}^{*}$}
\author[M.\ Carr]{MacKenzie Carr}
\author[K.\ Chaudhary]{Ketan Chaudhary}
\author[T.G.\ Marbach]{Trent G.\ Marbach}
\author[T.\ Mishura]{Teddy Mishura}

\address[A1,A2,A3]{Toronto Metropolitan University, Toronto, Canada}
\email[A1]{abonato@torontomu.ca}
\email[A2]{mackenzie.carr@torontomu.ca}
\email[A3]{ketan.chaudhary@torontomu.ca}

\address[A4]{Acadia University, Wolfville, Canada}
\email[A4]{trent.marbach@acadiau.ca}

\address[A5]{University of Delaware, Newark, USA}
\email[A5]{tmishura@udel.edu}

\thanks{${}^{*}$Supported by an NSERC Discovery Grant.}


\DeclareMathOperator{\ILMT}{ILMT}
\DeclareMathOperator{\ILTO}{ILTO}

\newcommand{\iltovec}[1]{\begin{bmatrix} L_{#1} \\ w_{#1} \\ \ell_{#1} \\ m_{#1} \\ a_{#1} \\ b_{#1} \\ e_{#1} \\ n_{#1} \end{bmatrix}}

\begin{document}

\begin{abstract}
Transitivity in social and other complex networks describes the tendency for two nodes with a common neighbor to become adjacent. We propose a new model for complex networks based on transitivity, called the Iterated Local Transitivity Model for Oriented Graphs (ILTO). The model generalizes the Iterated Local Model for Tournaments (ILMT). As in ILMT, oriented graphs are generated iteratively by cloning nodes and their adjacencies, with arcs between clones oriented either as the corresponding arcs between their parents or in the reverse direction. This process produces dense oriented graphs and, under suitable conditions, graphs with high connectivity and small diameter.

We analyze the universality properties of ILTO, showing that, for generating sequences with infinite support, the base oriented graph determines which oriented graphs appear as induced subgraphs, unlike in ILMT. We study motifs in ILTO and show that when the base graph is a tournament and the generating sequence has infinite support, the model generates quasirandom tournament sequences. Finally, we study domination numbers, the cop number, and the dichromatic number. We bound the cop number in terms of the independence number of the base graph and show that, for each fixed base graph, the dichromatic number is $\Theta(2^t)$ after $2t$ consecutive 0-steps.
\end{abstract}

\maketitle

\section{Introduction}

Complex networks are pervasive in the real world. Their observed properties include the small-world effect, scale-free degree distributions, and clustering into dense communities. Although these networks are often globally sparse, their dense communities motivate the study of densifying models. See the book \cite{bbook} for a survey of complex networks. 

We consider a densifying network model that generalizes the Iterated Local Model for Tournaments (ILMT) introduced and studied in \cite{ILMT}. Oriented graphs are directed graphs in which each pair of nodes shares at most one directed edge (or arc). Oriented graphs can represent social networks in which there are asymmetric interactions, without requiring every pair
of users to interact. For example, in an active sub-thread on the discussion site Reddit, in which users comment on one another's posts, an arc $(u,v)$ is assigned if user $u$ responds more frequently to the posts of user $v$ than $v$ does to $u$. However, if user $w$ has no interactions with user $v$ on the platform, then we do not assign an arc between $v$ and $w$; further, no arc is assigned when the two response frequencies are equal. Dense communities may form around popular or controversial users on the platform. Similar asymmetric communities also occur on other social media platforms, such as X, Instagram, and TikTok, or in sports, where arcs represent one player or team winning over another, with some players or teams possibly not in direct competition, such as in the FIFA World Cup. 

\emph{Balance theory} describes mechanisms for completing triads in social and other complex networks \cite{easl,hei}. One such
mechanism is \emph{transitivity}: if $x$ is a friend of $y$ and $y$ is a friend of $z$, then $x$ is a friend of $z$; see,
for example, \cite{scott}. Directed models of trust propagation were studied in \cite{guha}. A related framework, \emph{status theory}, interprets a positive arc as indicating that its creator views its recipient as having higher status. Evidence for this framework was found in directed networks derived from Epinions, Slashdot, and Wikipedia \cite{lesk}. 

Many models, including preferential attachment models \cite{ba,bol}, simulate properties of complex networks. The \emph{Iterated Local Transitivity (ILT)} model, introduced in \cite{ILT} and further studied in \cite{ILM,ILAT}, simulates structural properties in complex networks that arise from transitivity. Transitivity gives rise to the notion of \emph{cloning}, in which a new node $x$ is adjacent to an existing node $y$ and all the neighbors of $y$. As shown in \cite{ILT}, the model-generated graphs densify over time and exhibit poor spectral expansion. In addition, the ILT model generates graphs with the small-world property, exhibiting low diameter and a high clustering coefficient relative to random graphs with the same number of nodes and the same expected average degree. A directed analog of the ILT model was introduced in \cite{ILDT}, and a tournament-based model was first considered in \cite{ILTT}. Other variants of the ILT model include those found in \cite{ILH,FM,IGM,mscthesis,itind}. 

Tournament models require an interaction between every pair of nodes. Allowing nonadjacent pairs gives a broader class of networks, but also raises new questions about how the initial pattern of interactions constrains the graphs generated by the model. We extend the ILMT model to oriented graphs. The \emph{Iterated Local Transitivity Model for Oriented Graphs} (or \emph{ILTO}) is defined as follows. Fix an oriented graph $G_0$, called the \emph{base graph}, and an infinite binary sequence $s$, called the \emph{generating sequence}. We define a sequence of oriented graphs $G_t$ (where $t$ is referred to as the \emph{time-step}) as follows. 

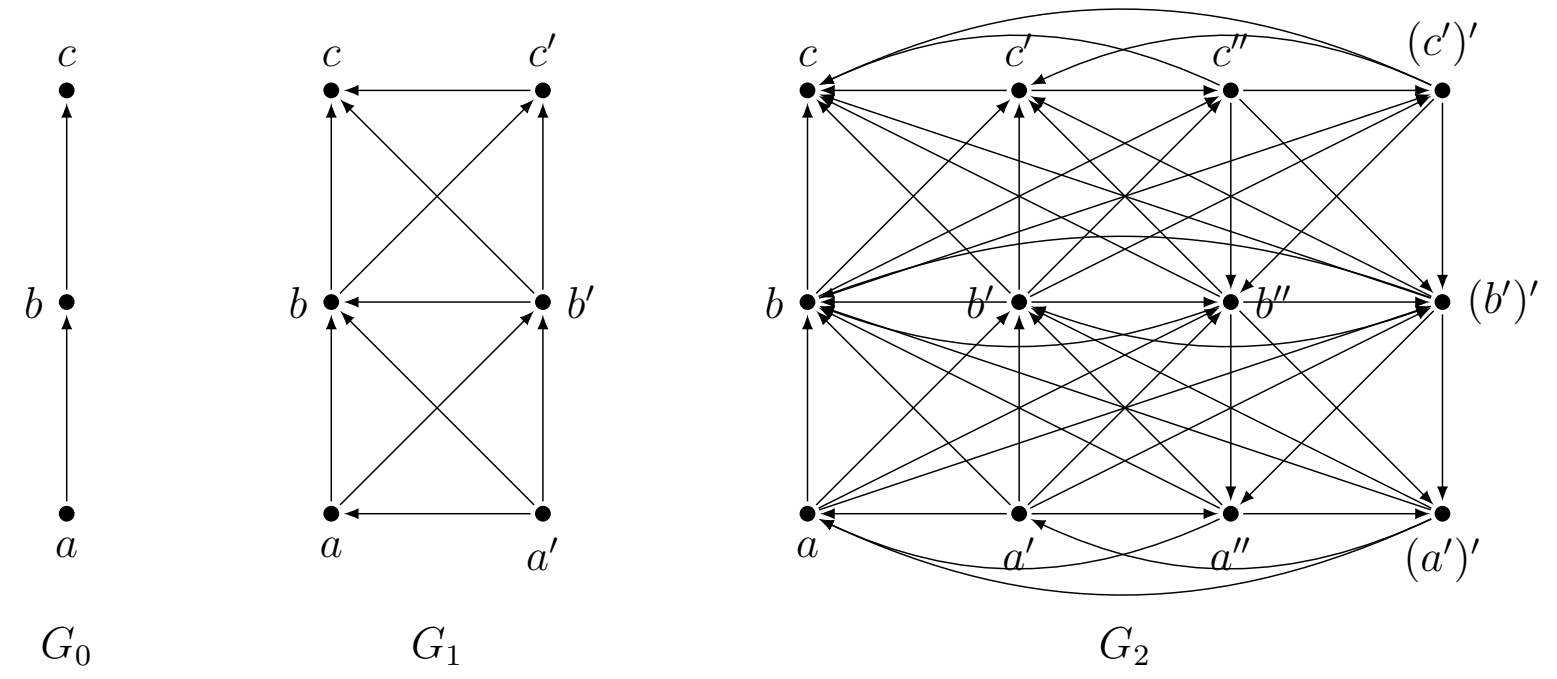
\begin{figure}[ht]
    \centering
    \scalebox{0.7}{\begin{tikzpicture}[
        vertex/.style={circle, fill=black, inner sep=1.5pt},
        edge/.style={->, >={latex[length=2.5mm, width=1.5mm]}, shorten >=1pt, shorten <=1pt}
    ]

    \begin{scope}
        \node[vertex, label=below:$a$] (a) at (0,0) {};
        \node[vertex, label=left:$b$] (b) at (0,2) {};
        \node[vertex, label=above:$c$] (c) at (0,4) {};

        \path[edge] (a) edge (b);
        \path[edge] (b) edge (c);

        \node at (0,-1.25) {$G_0$};
    \end{scope}

    \begin{scope}[xshift=2.5cm]
        \node[vertex, label=below:$a$] (a) at (0,0) {};
        \node[vertex, label=left:$b$] (b) at (0,2) {};
        \node[vertex, label=above:$c$] (c) at (0,4) {};
        \node[vertex, label=below:$a'$] (ap) at (2,0) {};
        \node[vertex, label=right:$b'$] (bp) at (2,2) {};
        \node[vertex, label=above:$c'$] (cp) at (2,4) {};

        \path[edge] (ap) edge (a);
        \path[edge] (bp) edge (b);
        \path[edge] (cp) edge (c);
        \path[edge] (a) edge (b);
        \path[edge] (ap) edge (b);
        \path[edge] (a) edge (bp);
        \path[edge] (ap) edge (bp);
        \path[edge] (b) edge (c);
        \path[edge] (bp) edge (c);
        \path[edge] (b) edge (cp);
        \path[edge] (bp) edge (cp);

        \node at (1,-1.25) {$G_1$};
    \end{scope}

    \begin{scope}[xshift=7cm]
        \node[vertex, label=below:$a$] (a) at (0,0) {};
        \node[vertex, label=left:$b$] (b) at (0,2) {};
        \node[vertex, label=above:$c$] (c) at (0,4) {};
        \node[vertex, label=below:$a'$] (ap) at (2,0) {};
        \node[vertex, label=left:$b'$] (bp) at (2,2) {};
        \node[vertex, label=above:$c'$] (cp) at (2,4) {};
        \node[vertex, label=below:$a''$] (app) at (4,0) {};
        \node[vertex, label=right:$b''$] (bpp) at (4,2) {};
        \node[vertex, label=above:$c''$] (cpp) at (4,4) {};
        \node[vertex, label=below:$(a')'$] (apprime) at (6,0) {};
        \node[vertex, label=right:$(b')'$] (bpprime) at (6,2) {};
        \node[vertex, label=above:$(c')'$] (cpprime) at (6,4) {};

        \path[edge] (ap) edge (a);
        \path[edge] (bp) edge (b);
        \path[edge] (cp) edge (c);
        \path[edge] (a) edge (b);
        \path[edge] (ap) edge (b);
        \path[edge] (a) edge (bp);
        \path[edge] (ap) edge (bp);
        \path[edge] (b) edge (c);
        \path[edge] (bp) edge (c);
        \path[edge] (b) edge (cp);
        \path[edge] (bp) edge (cp);

        \path[edge] (app) edge (apprime);
        \path[edge] (bpp) edge (bpprime);
        \path[edge] (cpp) edge (cpprime);
        \path[edge] (bpp) edge (app);
        \path[edge] (bpp) edge (apprime);
        \path[edge] (bpprime) edge (app);
        \path[edge] (bpprime) edge (apprime);
        \path[edge] (cpp) edge (bpp);
        \path[edge] (cpp) edge (bpprime);
        \path[edge] (cpprime) edge (bpp);
        \path[edge] (cpprime) edge (bpprime);

        \path[edge] (app) edge[bend left=25] (a);
        \path[edge] (bpp) edge[bend left=20] (b);
        \path[edge] (apprime) edge[bend left=25] (ap);
        \path[edge] (bpprime) edge[bend left=20] (bp);
        \path[edge] (cpp) edge[bend right=25] (c);
        \path[edge] (cpprime) edge[bend right=25] (cp);

        \path[edge] (apprime) edge[bend left=25] (a);
        \path[edge] (bpprime) edge[bend right=20] (b);
        \path[edge] (cpprime) edge[bend right=25] (c);
        \path[edge] (app) edge (b);
        \path[edge] (apprime) edge (b);
        \path[edge] (app) edge (bp);
        \path[edge] (apprime) edge (bp);
        \path[edge] (bpp) edge (c);
        \path[edge] (bpprime) edge (c);
        \path[edge] (bpp) edge (cp);
        \path[edge] (bpprime) edge (cp);

        \path[edge] (ap) edge (app);
        \path[edge] (bp) edge (bpp);
        \path[edge] (cp) edge (cpp);
        \path[edge] (a) edge (bpp);
        \path[edge] (ap) edge (bpp);
        \path[edge] (a) edge (bpprime);
        \path[edge] (ap) edge (bpprime);
        \path[edge] (b) edge (cpp);
        \path[edge] (bp) edge (cpp);
        \path[edge] (b) edge (cpprime);
        \path[edge] (bp) edge (cpprime);

        \node at (3,-1.25) {$G_2$};
    \end{scope}
    \end{tikzpicture}}
    \caption{The ILTO graphs $G_0$, $G_1$, and $G_2$, where $G_0$ is a directed path on three nodes and the generating sequence satisfies $s(1)=1$ and $s(2)=0$. In $G_2$, $x''$ denotes the clone of $x$ while $(x')'$ denotes the clone of $x'$, for $x=a,b,c$.}
    \label{fig:ILTO}
\end{figure}

For each time-step $t\geq 1$, we obtain the oriented graph $G_t$ from $G_{t-1}$ by the following steps. 
\begin{enumerate}
    \item For each node $x$ of $G_{t-1}$, add a node $x'$. 
    \item For each arc $(u,v)$ in $G_{t-1}$, add the arcs $(u,v')$ and $(u',v)$. 
    \item For each node $x$ in $G_{t-1}$, add the arc $(x',x)$. 
    \item For each arc $(u,v)$ in $G_{t-1}$, add the arc $(u',v')$ if $s(t)=1$ or add the arc $(v',u')$ if $s(t)=0$. 
\end{enumerate}
We denote $G_t$ more precisely by $\ILTO_{t,s}(G_0)$. For brevity, we often omit explicit references to $G_0$ and the generating sequence $s$. We refer to $x'$ as the \emph{clone} of $x$ in $G_t$ and $x$ as the parent of $x'$. The oriented graphs $G_t$ are called \emph{ILTO graphs}. If $s(t)=1$, then the process of obtaining $G_t$ from $G_{t-1}$ is referred to as a \emph{1-step}. A \emph{0-step} is defined analogously. See Figure~\ref{fig:ILTO} for examples of ILTO graphs. While a 1-step preserves the orientation between descendants of distinct nodes, a 0-step can create directed cycles even from an acyclic base graph: each arc $(u,v)$ gives the directed cycle $u,v',u',u$.

The subgraph of $G_t$ induced by the set of clone nodes is the \emph{reverse} of $G_{t-1}$ if $s(t)=0$ and is isomorphic to $G_{t-1}$ otherwise. We define the \emph{support} of $s$ to be $\{t\geq 1:s(t)=0\}$. If this set is infinite, then we say
that $s$ has \emph{infinite support}.

The paper is organized as follows. In Section~2, we prove that ILTO graphs densify over time and, under suitable conditions, have high connectivity and small diameter. Under infinite support, every finite tournament eventually appears as an induced subgraph in either model. For general oriented graphs, however, induced containment in ILTO depends on the underlying base graph. For a given oriented graph $H$, we characterize the base graphs $G_0$ for which $\ILTO_{t,s}(G_0)$ contains an induced copy of $H$ for large $t$ when $s$ has infinite support; see Theorem~\ref{thm:ILTOuni}. In Section~3, we give precise counts of tournaments of order at most 4 in an ILTO graph and prove that, when the base graph is a tournament and the generating sequence has infinite support, the resulting graphs form a quasirandom sequence; see Theorem~\ref{thm:0stepquasi}.

In Section~4, we consider graph-theoretical properties of ILTO graphs, including their domination number, cop number, and dichromatic number. We construct a tournament $G_0$ with in-domination number 2 such that $G_1$ has in-domination number 3 when $s(1)=0$. We characterize the in-dominating sets $S$ in the oriented graph $G_0$ for which $S' = \{x' : x\in S\}$ is in-dominating for $G_1$ when $s(1)=0$; see Lemma~\ref{lem:indom}. We also show that the in-domination number is bounded in terms of the base graph, independently of the generating sequence and the number of time-steps; see Theorem~\ref{thm:bounded indomination}. We prove that, for each fixed base graph $G_0$, the dichromatic number after $2t$ consecutive 0-steps is $\Theta(2^t)$; see Theorems~\ref{thm:upper coloring} and~\ref{thm:lower coloring}. Finally, in Section~5, we discuss directions for future research. 

Throughout the paper, we consider finite, simple, and connected oriented graphs. An oriented graph $G$ is \emph{connected} if its underlying undirected graph, which we denote by $G^*$, is connected. An oriented graph is \emph{strongly connected} if for each pair of nodes $x$ and $y$, there exist directed paths connecting $x$ to $y$ and $y$ to $x$. The \emph{in-neighbors} and \emph{out-neighbors} of a node $v$ are denoted by $N^-(v)$ and $N^+(v)$, respectively. Similarly, the \emph{in-degree} and \emph{out-degree} of $v$ are denoted by $d^-(v)$ and $d^+(v)$, respectively. An oriented graph is \emph{acyclic} if it does not contain a directed cycle as a subgraph. A \emph{sink} is a node with out-degree 0, and a \emph{source} is a node with in-degree 0. For background on graphs and digraphs, the reader is directed to \cite{west} and \cite{bang-jensen-gutin-2009}, respectively. For background on social and complex networks, see \cite{bbook,at}.

\section{Structural Properties and Universality}

In this section, we show that ILTO graphs satisfy several properties observed in real-world complex networks, including densification and, under suitable conditions, high connectivity and small diameter. For a given base graph $G_0$, we denote the number of nodes in $G_t = \ILTO_{t,s}(G_0)$ by $n_t$ and the number of arcs by $e_t$. We begin with a simple formula for these parameters in terms of the base graph's order and size. 

\begin{theorem}
\label{prop:n and e}
    Let $G_0$ be an oriented graph with $n_0$ nodes and $e_0$ arcs. We then have that
    \begin{enumerate}
        \item $n_t = 2^tn_0$, and 
        \item $e_t = 4^te_0 + 2^{t-1}(2^t-1)n_0.$
    \end{enumerate}
\end{theorem}

\begin{proof}
    Item (1) follows directly from the ILTO model definition. For item (2), when performing either a 0-step or a 1-step, each arc $(u,v)$ in $G_{t-1}$ produces four arcs in $G_t$: $(u,v)$, $(u',v)$, $(u,v')$ and exactly one of $(u',v')$ or $(v',u')$. Additionally, for each of the $n_{t-1}$ nodes $u$ in $G_{t-1}$, we have the arc $(u',u)$ in $G_t$. Thus, $e_t = 4e_{t-1}+n_{t-1}$. 

    Using $n_{t-1}=2^{t-1}n_0$ and iterating the recurrence, we obtain
$$e_t=4^te_0+n_0\sum_{i=0}^{t-1}4^{t-1-i}2^i
=4^te_0+2^{t-1}(2^t-1)n_0,$$
as required.
\end{proof}

We now show a densification power law for ILTO graphs. For positive real-valued sequences $f_t$ and $g_t$, we write
$f_t\sim g_t$ if $f_t/g_t$ tends to 1 as $t\to\infty$.

\begin{theorem}
Let $G_0$ be an oriented graph with $n_0$ nodes and $e_0$ arcs, and let $G_t=\ILTO_{t,s}(G_0)$. As $t\to\infty$, we have that
$$\frac{e_t}{n_t}\sim 2^t\left(\frac{2e_0+n_0}{2n_0}\right).$$
In particular, $e_t\sim Cn_t^2$, where $C=\frac{2e_0+n_0}{2n_0^2}$.
\end{theorem}

\begin{proof}
    From above, we have that $e_t = 4^te_0 + 2^{t-1}(2^t-1)n_0 = 2^{2t-1}\left(2e_0\right) + 2^{2t-1}n_0-2^{t-1}n_0$. We then have that
 \begin{eqnarray*}
\frac{e_t}{n_t} &\sim& \frac{2^{2t-1}(2e_0)+2^{2t-1}n_0}{2^tn_0}\\
&=& 2^t\left(\frac{2e_0+n_0}{2n_0}\right).
\end{eqnarray*}
Since $n_t=2^tn_0$, we also have $e_t\sim\frac{2e_0+n_0}{2n_0^2}n_t^2$, proving the second assertion.
\end{proof} 

Next, we compare the number of arcs in an ILTO graph to that of a tournament of the same order. We show that this ratio is nondecreasing over time and tends to 1 if and only if the base graph is a tournament.

\begin{theorem}
Let $G_0$ be an oriented graph and $G_t=\ILTO_{t,s}(G_0)$. Let $e^*(t)=\binom{2^tn_0}{2}$. For all $t\geq1$ such that $n_{t-1}\geq2$, we have that $$\frac{e_t}{e^*(t)}\geq\frac{e_{t-1}}{e^*(t-1)}.$$
Further, $\lim_{t\to\infty}\frac{e_t}{e^*(t)}=1$ if and only if $G_0$ is a tournament.
\end{theorem}

\begin{proof}
    If $n_0=1$, then every $G_t$ is a tournament, and $e_t/e^*(t)=1$ for all $t\geq 1$. Thus, both conclusions hold. We may therefore assume that $n_0\geq 2$. A straightforward calculation using Theorem~\ref{prop:n and e} shows that 
    $$\frac{e_t}{e^*(t)}-\frac{e_{t-1}}{e^*(t-1)} = \frac{2^{2t-2}\left(n_0^2-n_0-2e_0\right)}{2^{t-1}n_0(2^tn_0-1)(2^{t-1}n_0-1)}.$$ 
    Since $e_0 \leq \binom{n_0}{2} = \frac{n_0^2-n_0}{2}$, we have that $$\frac{e_t}{e^*(t)}-\frac{e_{t-1}}{e^*(t-1)}\geq 0.$$

    Thus, $$\lim_{t\to \infty} \frac{e_t}{e^*(t)} = \lim_{t\to \infty} \frac{2^{t+1}e_0 + 2^t n_0 -n_0}{2^tn_0^2-n_0} = \frac{2e_0+n_0}{n_0^2}.$$
    We have that $\frac{2e_0+n_0}{n_0^2}=1$ if and only if $e_0 = \frac{n_0^2-n_0}{2} = \binom{n_0}{2}$. This is true exactly when $G_0$ is a tournament. 
\end{proof}

We next consider the connectivity of ILTO graphs. For an oriented graph $G$ with at least two nodes, the \emph{connectivity} $\kappa(G)$ is the minimum cardinality of a set of nodes $S\subseteq V(G)$ such that $G\setminus S$ is not strongly connected. We set $\kappa(K_1)=0$. An oriented graph $G$ is called \emph{$k$-connected} if $\kappa(G)\geq k$. We show that if $G_0$ has a connected underlying graph and no sink node, then a 0-step produces a strongly connected graph, even if $G_0$ is not strongly connected. 

\begin{theorem}\label{thm:zero_step_strongly_connects}
If $G_0$ is an oriented graph with no sink node and $s(1)=0$, then $G_1=\ILTO_{1,s}(G_0)$ is strongly connected.
\end{theorem}

\begin{proof}
For each $u\in V(G_0)$, choose an out-neighbor $v$, which exists since $G_0$ has no sink. The arcs $(u,v')$, $(v',u')$, and $(u',u)$ form a directed 3-cycle in $G_1$. Thus, $u$ and its clone $u'$ belong to the same strongly connected component of $G_1$.

For every arc $(u,v)$ in $G_0$, the arcs $(u,v)$ and $(v',u')$ belong to $G_1$. Since each parent and its clone belong to the same strongly connected component, these two arcs imply that $u,u',v,v'$ all belong to the same strongly connected component. As $G_0^*$ is connected, all nodes of $G_1$ therefore belong to a single strongly connected component.
\end{proof}

By the following theorem, once an ILTO graph with at least two nodes is strongly connected, its connectivity grows without bound as the number of time-steps increases. We note that this is a generalization of \cite[Theorem~2]{ILMT}, whose full proof was not included in \cite{ILMT}. 

\begin{theorem}\label{thm: kconnectivity}
If $G_0$ is a $k$-connected oriented graph, for some $k\geq 1$, and $G_1 = \mathrm{ILTO}_{1,s}(G_0)$ for some generating sequence $s$, then $G_1$ is $2k$-connected.
\end{theorem}

\begin{proof}
Since $G_0$ is $k$-connected, for each $v\in V(G_0)$, we have $d_{G_0}^-(v),d_{G_0}^+(v)\geq k$. It follows from the ILTO construction that each node of $G_1$ has at least $k$ in-neighbors and at least $k$ out-neighbors in each of $V(G_0)$ and $V(G_0)' = V(G_1)\setminus V(G_0)$. Let $S\subseteq V(G_1)$ with $|S|\leq 2k-1$. By the Pigeonhole Principle, either $|S\cap V(G_0)|\leq k-1$ or $|S\cap V(G_0)'|\leq k-1$.

In the former case, each node in $V(G_1)\setminus S$ has at least one in-neighbor and one out-neighbor in $V(G_0)\setminus S$. Since $G_0$ is $k$-connected, the subgraph induced by $V(G_0)\setminus S$ is strongly connected. Every clone node $v'\in V(G_0)'\setminus S$ inherits all of its parent's in- and out-neighbors in $V(G_0)$ and also has its parent as an out-neighbor. Since at most $k-1$ parents were deleted, $v'$ has both an in-neighbor and an out-neighbor in $V(G_0)\setminus S$. Therefore, $G_1\setminus S$ is strongly connected. The latter case follows similarly, since the subgraph induced by $V(G_0)'$ is isomorphic to either $G_0$ or its reverse and is therefore also $k$-connected. Thus, $G_1$ is $2k$-connected.
\end{proof}

We next show that after two 0-steps, the diameter of an ILTO graph is at most two more than the diameter of the underlying base graph, provided the base oriented graph has no sink. In the proof, we use the notation $u \rightarrow v$ to mean that $(u,v)$ is an arc.

\begin{theorem}\label{thm:zeroupper}
Let $G_0$ be an oriented graph with no sink node, and let $s$ be a generating sequence containing at least two 0-entries. If $r$ is the index of the second 0-entry, then for all $t\geq r$, $$\mathrm{diam}(\ILTO_{t,s}(G_0))\leq \mathrm{diam}(G_0^*)+2.$$
\end{theorem}

\begin{proof}
Write $G_t=\ILTO_{t,s}(G_0)$ and let $D=\mathrm{diam}(G_0^*)$. Cloning preserves the diameter of the underlying graph and creates no sinks. After the first 0-step, every node lies on a directed 3-cycle: for each node $u$, an out-neighbor $v$ gives the cycle $u\to v'\to u'\to u$. This property persists, since replacing any one node of a directed 3-cycle by its clone gives another directed 3-cycle.

First suppose that $s(1)=s(2)=0$. For each $u\in V(G_0)$, let $S_u=\{u',u'',(u')'\}$, where $u'$ is created at the first step, and $u''$ and $(u')'$ are the clones of $u$ and $u'$ at the second step. Each $S_u$ induces a directed 3-cycle, and all its nodes point to $u$. If $(u,v)$ is an arc of $G_0$, then the arcs from $S_u$ to $S_v$ are $(u',v'')$, $(u'',v')$, and $((u')',(v')')$, with all other cross arcs directed oppositely. Thus, if $u,v$ are adjacent in $G_0^*$, then every node of either set has an in-neighbor and an out-neighbor in the other.

It follows that a path $u_0,u_1,\ldots,u_k$ in $G_0^*$ gives a directed path $z_0,z_1,\ldots,z_k$ with $z_i\in S_{u_i}$, where either endpoint may be prescribed. If $x\in S_u$ and $y\in S_v$, then choose a shortest underlying path from $u$ to $v$ and a corresponding directed path ending at $y$. Its first node can be reached from $x$ in at most two arcs within $S_u$, giving distance at most $D+2$. If instead $y=v$, then choose the directed path to start at $x$ and append its last node's arc to $v$, giving distance at most $D+1$. Finally, if $x=u\in V(G_0)$, then choose an out-neighbor $w$ of $u$. Since $u$ points to every node of $S_w$, a shortest underlying path from $w$ to $v$ gives distance at most $D+1$ to any prescribed node of $S_v$, and at most $D+2$ to $v$. Hence, $\mathrm{diam}(G_2)\leq D+2$.

Now let $q<r$ be the indices of the first two 0-steps, and fix distinct nodes $x,y\in V(G_r)$. For a node of $G_r$, its \emph{ancestor in $G_q$} is obtained by repeatedly replacing a clone by its parent until reaching a node of $G_q$. A node already in $G_q$ is its own ancestor. We call each node a \emph{descendant} of its ancestor. If the ancestors of $x$ and $y$ in $G_q$ are distinct, then choose one descendant in $G_{r-1}$ of each node of $G_q$, choosing $x$ and $y$ themselves when they belong to $G_{r-1}$, and their parents otherwise. Since all intervening steps are 1-steps, the chosen nodes induce a copy of $G_q$. This copy and its clones in $G_r$ contain $x,y$ and induce the graph obtained from $G_{q-1}$ by two consecutive 0-steps. The preceding argument applies, since $G_{q-1}$ has no sink and $\mathrm{diam}(G_{q-1}^*)=D$.

If both ancestors are the same node $u$, then choose a directed 3-cycle $u\to v\to w\to u$ in $G_q$. The construction preserves the arcs $x\to v$, $v\to w$, and $w\to y$, giving a path of length three. Therefore, $\mathrm{diam}(G_r)\leq D+2$.

Finally, once every node lies on a directed 3-cycle, any diameter bound of at least three is preserved by a time-step. Paths of length at least two remain directed paths when either endpoint or both endpoints are replaced by their clones. An arc requires at most two arcs after these replacements, using $a'\to a\to b'$ when necessary. A clone reaches its parent directly, while a cycle $a\to b\to c\to a$ gives the path $a\to b\to c\to a'$ in the opposite direction. Since $D+2\geq3$, the bound holds for every $t\geq r$.
\end{proof}

When the base graph $G_0$ is a tournament, the following result of \cite{ILMT} requires only nonempty support.

\begin{corollary}
    \cite[Theorem~1]{ILMT} If the tournament $G_0$ has no sink node and $s$ has nonempty support, then for sufficiently large $t$, the diameter of $\ILMT_{t,s}(G_0)$ is at most 3.
\end{corollary}

We next investigate which oriented graphs appear as subgraphs of ILTO graphs, for a given base graph. In Theorem~4 of \cite{ILMT}, we show the following universality property for the ILMT model.

\begin{theorem} [\cite{ILMT}]
Let $G_0$ be a tournament of order $n_0$ and let $n\leq n_0$. For a generating sequence $s$ with infinite support, let $r$ be the smallest integer such that $$\left|\{1 \leq t \leq r : s(t)=0\}\right| = n.$$ Hence, every tournament on $n$ nodes is isomorphic to a subtournament of $\mathrm{ILMT}_{r,s}(G_0)$.
\label{lemma:0steps}
\end{theorem}

Unlike the tournament case, not every oriented graph appears as an induced subgraph of $\ILTO_{t,s}(G_0)$ for some $t$, even when $s$
has infinite support. A simple example of this is if $G_0$ is a tournament on three nodes and we consider any orientation of the path on three nodes. For all $t\geq 1$, $\ILTO_{t,s}(G_0)$ is a tournament, and thus any three nodes induce either a directed or transitive 3-cycle, not an orientation of a path. However, for a given oriented graph $H$, we can identify properties of the base graph $G_0$ that imply that $H$ is an induced subgraph of $\ILTO_{t,s}(G_0)$ for some value of $t$. 

For $t>0$ and $G_t = \ILTO_{t,s}(G_0)$, define the function $p_t:V(G_t) \to V(G_{t-1})$ given by $$p_t(v) = \begin{cases}
    v & \text{if } v\in V(G_{t-1}) \\
    u & \text{if } v=u'\text{ and } u\in V(G_{t-1}).
\end{cases}$$ Let $p: V(G_t) \to V(G_0)$ be given by $p(v) = (p_1 \circ p_{2}\circ \dots \circ p_t)(v)$ and call this the \emph{projection of $v$ onto $G_0$}. The following lemma characterizes adjacency in $G_t$ in terms of the projection onto $G_0$.

\begin{lemma}
    \label{lem:projection}
    Let $G_0$ be an oriented graph, $G_t = \ILTO_{t,s}(G_0)$ and $p: V(G_t) \to V(G_0)$ be defined as above. For $u,v\in V(G_t)$ with $u\neq v$, there is an arc between nodes $u$ and $v$ if and only if $p(u)=p(v)$ or there is an arc in $G_0$ between nodes $p(u)$ and $p(v)$. 
\end{lemma}
\begin{proof}
    Let $t=1$, and note that $p(v) = p_1(v)$ for all $v\in V(G_1)$. If $p(u)=p(v)$, then either $u=v'$ or $v=u'$. In the former case, we have the arc $(u,v)$, and in the latter case, we have the arc $(v,u)$ in $G_1$. 

    Now assume, without loss of generality, that $(p(u),p(v))$ is an arc in $G_0$. In $G_1$, we have the arcs $(p(u),p(v))$, $(p(u)',p(v))$, $(p(u),p(v)')$ and exactly one of $(p(u)',p(v)')$ or $(p(v)',p(u)')$. By definition, we have $u\in \{p(u),p(u)'\}$ and $v\in \{p(v),p(v)'\}$. Thus, there is an arc between nodes $u$ and $v$ in $G_1$. 

    Finally, if $p(u)\neq p(v)$ and there is no arc between $p(u)$ and $p(v)$ in $G_0$, then in the construction of $G_1$, no arc is added between $p(u)$ and $p(v)$, between $p(u)'$ and $p(v)$, between $p(u)$ and $p(v)'$, or between $p(u)'$ and $p(v)'$. Thus, there is no arc between $u$ and $v$ in $G_1$. 

    The result now follows by induction on $t$, by setting $G_0' = G_{t-1}$ and $G'_1=G_t$ and repeating the above argument for $G'_0$ and $G'_1$. 
\end{proof}

Next, we define a subgraph of the underlying graph of $G$ to determine exactly which ILTO graphs contain $G$ as an induced subgraph. Given an oriented graph $G$, from the underlying graph $G^*$, we construct the (undirected) graph $N^*_G$ as follows. For $u\in V(G^*)$, let $N_{G^*}[u]$ denote the closed neighborhood of $u$, consisting of $u$ and its neighbors in $G^*$. For each pair $u,v\in V(G^*)$, define the equivalence relation $u\sim v$ if $N_{G^*}[u] = N_{G^*}[v]$. Let $V_1,V_2,\dots,V_m$ be the equivalence classes of $V(G^*)$ with respect to $\sim$. Note that, in $G^*$, the set of edges $E[V_i,V_j]$ between equivalence classes $V_i$ and $V_j$ is either empty or complete, for all $i,j\in \{1,2,\dots,m\}$, $i\neq j$. Now, let $N^*_G$ be the graph on the node set $\{V_1,V_2,\dots,V_m\}$ where $V_i$ and $V_j$ are adjacent in $N^*_G$ if and only if $E[V_i,V_j]$ is complete. Equivalently, $N^*_G$ can be obtained from $G^*$ by deleting all but one node from $V_i$, for $i=1,2,\dots,m$. 

We now use this construction to characterize the base graphs whose ILTO sequences eventually contain $G$ as an induced subgraph, provided the generating sequence has infinite support.

\begin{theorem}\label{thm:ILTOuni}
    Let $H$ be an oriented graph, $s$ a binary sequence with infinite support, and $N^*_H$ the undirected graph constructed as above. Given an oriented graph $G_0$, $N^*_H$ is an induced subgraph of $G_0^*$ if and only if there exists $t_0=t_0(H,s)$ such that for all $t\geq t_0$, $H$ is an induced subgraph of $\ILTO_{t,s}(G_0)$. 
\end{theorem}

\begin{proof}
Suppose that $N^*_H$ is an induced subgraph of $G_0^*$. Let $V_1,V_2,\dots,V_m$ be the equivalence classes used in the construction of $N^*_H$ and label the nodes of an induced copy of $N^*_H$ in $G_0^*$ by their corresponding equivalence classes. Let $q=\left\lceil\log_2\left(\max_{1\leq i\leq m}|V_i|\right)\right\rceil$. We show next that $H^*$ is an induced subgraph of $G_q^*$, where $G_q=\ILTO_{q,s}(G_0)$. For $v\in V_i$, we have
$$N_{H^*}[v]=\left(\bigcup_{V_iV_j\in E(N^*_H)}V_j\right)\cup V_i.$$

For $v\in V(G_0)$, let $S_v=\{x\in V(G_q):p(x)=v\}$, where $p:V(G_q)\to V(G_0)$ is the projection onto $G_0$, taking $p$ to be the identity if $q=0$. Observe that $|S_v|=2^q\geq\max_{1\leq i\leq m}|V_i|$ for every $v\in V(G_0)$. For $i=1,2,\dots,m$, choose $S_i\subseteq S_{V_i}$ with $|S_i|=|V_i|$ and let $S=\bigcup_{i=1}^m S_i$.

For distinct $x,y\in S_i$, we have $p(x)=p(y)=V_i$. By Lemma~\ref{lem:projection}, there is an arc between $x$ and $y$. Thus, $S_i$ induces a clique in $G_q^*$. For $x\in S_i$ and $y\in S_j$ with $i\neq j$, Lemma~\ref{lem:projection} gives that $x$ and $y$ are adjacent in $G_q^*$ if and only if $V_i$ and $V_j$ are adjacent in $G_0^*$. Since the chosen nodes induce a copy of $N_H^*$, this holds if and only if $V_iV_j\in E(N_H^*)$. Therefore, the closed neighborhood of $x$ in $G_q^*[S]$ is
$$\left(\bigcup_{V_iV_j\in E(N_H^*)}S_j\right)\cup S_i.$$
Thus, $S$ induces a copy of $H^*$ in $G_q^*$. In $G_q$, it induces an orientation of $H^*$, with some arcs possibly reversed relative to $H$. We use additional 0-steps to correct these orientations.

Let $S^{(0)}=S$ and let $\phi_0:S^{(0)}\to V(H)$ be an isomorphism between the underlying graphs. Choose a node $x\in V(H)$, let $x_0=\phi_0^{-1}(x)$, and let $t_1=\min\{t>q:s(t)=0\}$. Let $A^{(0)}$ be the set of nodes $v\in S^{(0)}\setminus\{x_0\}$ for which either $(x_0,v)$ is an arc in $G_q$ and $(\phi_0(v),x)$ is an arc in $H$, or $(v,x_0)$ is an arc in $G_q$ and $(x,\phi_0(v))$ is an arc in $H$.

Form $S^{(1)}$ from $S^{(0)}$ by replacing the elements of $A^{(0)}\cup\{x_0\}$ by their clones in $G_{t_1}$. By the construction of the model, $S^{(1)}$ still induces a copy of $H^*$. Define $\phi_1:S^{(1)}\to V(H)$ by $\phi_1(v)=\phi_0(v)$ for nodes retained from $S^{(0)}$ and $\phi_1(v')=\phi_0(v)$ for nodes replaced by their clones. Cloning both endpoints reverses each incorrectly oriented arc incident with $x_0$, while cloning only $x_0$ preserves each correctly oriented arc incident with it. Nonadjacent pairs remain nonadjacent. Thus, all arcs incident with the node corresponding to $x$ now have the required orientations.

Repeat this process at successive 0-steps, each time choosing a node of $H$ that was not previously chosen and defining $S^{(i)}$ and $\phi_i$ analogously. No node corresponding to a previously chosen node of $H$ is cloned, since all its incident arcs already have the required orientations. These arcs therefore remain correct. Let $n_H=|V(H)|$. After $n_H$ additional 0-steps, the bijection $\phi_{n_H}:S^{(n_H)}\to V(H)$ preserves arcs and nonadjacent pairs. Therefore, $H$ is an induced subgraph of $G_{t_{n_H}}$ and remains so at every later time-step. Taking $t_0=t_{n_H}$ gives the required threshold, which depends only on $H$ and $s$.

Conversely, suppose that $H$ is an induced subgraph of $G_{t^*}=\ILTO_{t^*,s}(G_0)$ for some $t^*\geq0$, and let $V_H\subseteq V(G_{t^*})$ induce a copy of $H$. Retain one node from each equivalence class to obtain a copy of $N_H^*$ in $G_{t^*}^*$. Label these representatives $V_1,V_2,\dots,V_m$ by their corresponding classes. Let $p:V(G_{t^*})\to V(G_0)$ be the projection onto $G_0$, taking $p$ to be the identity if $t^*=0$.

Let $V_i,V_j$ be distinct representatives. By Lemma~\ref{lem:projection}, if $V_iV_j\notin E(N_H^*)$, then $p(V_i)\neq p(V_j)$ and $p(V_i)p(V_j)\notin E(G_0^*)$. Suppose instead that $V_iV_j\in E(N_H^*)$ and $p(V_i)=p(V_j)$. By Lemma~\ref{lem:projection}, a representative $V_\ell$ belongs to $N_{N_H^*}[V_i]$ precisely when $p(V_\ell)=p(V_i)$ or
$p(V_\ell)p(V_i)\in E(G_0^*)$. Since $p(V_i)=p(V_j)$, it follows that $N_{N_H^*}[V_i]=N_{N_H^*}[V_j]$. Each closed neighborhood in $H^*$ is the union of the classes represented in the corresponding closed neighborhood of $N_H^*$. Thus, nodes in the classes $V_i$ and $V_j$ have equal closed neighborhoods in $H^*$, contradicting that these classes are distinct.

Therefore, the projections of the representatives are distinct. By Lemma~\ref{lem:projection}, $p(V_i)p(V_j)\in E(G_0^*)$ if and only if $V_iV_j\in E(N_H^*)$. Hence, $\{p(V_i):1\leq i\leq m\}$ induces a copy of $N_H^*$ in $G_0^*$.
\end{proof}

For any oriented graph $H$, the undirected graph $N^*_H$ is an induced subgraph of $H^*$. Thus, we have the following corollary. 

\begin{corollary}
    Let $H$ be an oriented graph and $s$ a binary sequence with infinite support. For an oriented graph $G$, if $H^*$ is an induced subgraph of $G^*$, then there exists some $t_0 = t_0(H,s)$ such that for all $t\geq t_0$, $H$ is an induced subgraph of $\ILTO_{t,s}(G)$. 
\end{corollary}

We note that if $H$ is a tournament, then the graph $N^*_H$ is exactly the graph consisting of a single isolated node. Since $N^*_H$ is an induced subgraph of every nonempty undirected graph, $H$ is an induced subgraph of $\ILTO_{t,s}(G_0)$ for all sufficiently large $t$ and any oriented base graph $G_0$, provided that $s$ has infinite support.

\begin{corollary}\label{couni}
    Let $G_0$ be an oriented graph and $s$ a generating sequence with infinite support. We then have that every tournament on $n$ nodes is an induced subgraph of $G_t = \mathrm{ILTO}_{t,s}(G_0)$ for some $t \ge 0$.
\end{corollary}

Finally, because every oriented graph $H$ is a subgraph (not necessarily induced) of some tournament on $|V(H)|$ nodes, we have that $H$ appears as a subgraph of $\ILTO_{t,s}(G_0)$ for sufficiently large $t$. 

\begin{corollary}
    If $H$ is an oriented graph and $s$ is a binary sequence with infinite support, then there exists some $t_0 = t_0(H,s)$ such that for all $t\geq t_0$ and all oriented graphs $G$, $H$ is a subgraph of $\ILTO_{t,s}(G)$. 
\end{corollary}

\section{Motifs and Quasirandomness}

Having explored in the previous section when a particular oriented graph appears as an induced subgraph of ILTO graphs, we now study the counts of particular induced subgraphs in ILTO graphs. A \emph{motif} is an isomorphism type of a specified oriented graph. Motifs play an important role in the structure of complex networks, both in graphs and in digraphs; see \cite{Milo}. We provide precise counts of motifs of order at most 4 that are isomorphic to tournaments.

For an oriented graph $G_0$, let $G_t=\ILTO_{t,s}(G_0)$ for each $t\geq 0$. There are two tournaments on three nodes: the directed 3-cycle, $D_3$, and the transitive 3-cycle, $T_3$. Let $a_t$ denote the number of subsets of $V(G_t)$ of cardinality 3 that form a $D_3$ and $b_t$ denote the number of subsets of $V(G_t)$ that form a $T_3$. Similarly, there are four tournaments on four nodes: the linear order, $T_4$, the \emph{Winner} tournament, which is the unique non-transitive tournament with a node of out-degree 3, the \emph{Loser} tournament, which is the unique non-transitive tournament with a node of in-degree 3, and the \emph{Mixed} tournament, which is the unique tournament with out-degree sequence $(1,1,2,2)$. Let $L_t$ be the number of four-element subsets in $G_t$ that form a linear order, $w_t$ the number of four-element subsets in $G_t$ that form a Winner tournament, $\ell_t$ the number of four-element subsets in $G_t$ that form a Loser tournament, and $m_t$ the number of four-element subsets in $G_t$ that form a Mixed tournament.

\begin{theorem}\label{thm: motifs}
    Let $G_0$ be an oriented graph and $G_1 = \ILTO_{1,s}(G_0)$. If $s(1) = 1$, then
    $$\iltovec{1} = \begin{bmatrix}
    16 & 0 & 0 & 0 & 0 & 12 & 1 & 0 \\
    0 & 16 & 0 & 0 & 0 & 0 & 0 & 0\\
    0 & 0 & 16 & 0 & 0 & 0 & 0 & 0\\
    0 & 0 & 0 & 16 & 12 & 0 & 0 & 0\\
    0 & 0 & 0 & 0 & 8 & 0 & 0 & 0\\
    0 & 0 & 0 & 0 & 0 & 8 & 4 & 0\\
    0 & 0 & 0 & 0 & 0 & 0 & 4 & 1\\
    0 & 0 & 0 & 0 & 0 & 0 & 0 & 2
\end{bmatrix} \iltovec{0}.$$
    If $s(1)=0$, then
    $$\iltovec{1} = \begin{bmatrix}
    11 & 3 & 3 & 3 & 0 & 6 & 0 & 0 \\
    1 & 9 & 1 & 1 & 3 & 1 & 0 & 0\\
    1 & 1 & 9 & 1 & 0 & 2 & 1 & 0\\
    3 & 3 & 3 & 11 & 9 & 3 & 0 & 0\\
    0 & 0 & 0 & 0 & 5 & 1 & 1 & 0\\
    0 & 0 & 0 & 0 & 3 & 7 & 3 & 0\\
    0 & 0 & 0 & 0 & 0 & 0 & 4 & 1\\
    0 & 0 & 0 & 0 & 0 & 0 & 0 & 2
\end{bmatrix} \iltovec{0}.$$
\end{theorem}
\begin{proof}
Suppose that $s(1)=1$. A copy of $T_4$ that does not contain both a node $u$ and its clone $u'$ can be constructed from a copy of $T_4$ in $G_0$ by replacing any number of the four nodes by their clones, contributing $16L_0$ to $L_1$. Consider a copy of $T_3$ in $G_0$, with nodes labeled $u,v,w$. The subgraph of $G_1$ induced by $\{u,u',x,y\}$, where $x\in\{v,v'\}$ and $y\in\{w,w'\}$, forms a copy of $T_4$ in $G_1$. Similarly, $\{x,v,v',y\}$ where $x\in\{u,u'\}$ and $y\in\{w,w'\}$ forms a $T_4$ in $G_1$, as does $\{x,y,w,w'\}$ with $x\in\{u,u'\}$ and $y\in\{v,v'\}$. These contribute a total of $12b_0$ to $L_1$. Finally, consider an arc $(u,v)$ in $G_0$. The subgraph induced by $\{u,u',v,v'\}$ forms a copy of $T_4$ in $G_1$. Thus, $L_1=16L_0+12b_0+e_0$.

A copy of the Winner tournament in $G_1$ cannot contain both a node $u$ in $G_0$ and its clone $u'$. Therefore, any copy of the Winner tournament in $G_1$ corresponds to a copy of the Winner tournament in $G_0$ by replacing any number of the four nodes with their clones. Thus, $w_1 = 16w_0$. A similar argument shows that $\ell_1 = 16\ell_0$.

A copy of the Mixed tournament that does not contain both a parent node and its clone can be constructed from a copy of the Mixed tournament in $G_0$ by replacing any number of its nodes by their clones, contributing $16m_0$ to $m_1$. Consider a copy of $D_3$ in $G_0$ with nodes labeled $u,v,w$. The subgraph of $G_1$ induced by $\{u,u',x,y\}$, where $x\in\{v,v'\}$ and $y\in\{w,w'\}$, forms a copy of the Mixed tournament in $G_1$. Similarly, $\{x,v,v',y\}$ where $x\in\{u,u'\}$ and $y\in\{w,w'\}$ forms a Mixed tournament in $G_1$, as does $\{x,y,w,w'\}$ with $x\in\{u,u'\}$ and $y\in\{v,v'\}$. These contribute a total of $12a_0$ to $m_1$. Thus, $m_1=16m_0+12a_0$.

Note that no $D_3$ in $G_1$ can contain both a node $u$ in $G_0$ and its clone $u'$. Therefore, any copy of $D_3$ in $G_1$ corresponds to a copy of $D_3$ in $G_0$ by replacing any number of the three nodes with their clones. Thus, $a_1=8a_0$.

A copy of $T_3$ in $G_1$ that does not contain both a parent node and its clone can be constructed from a copy of $T_3$ in $G_0$ by replacing any number of the three nodes with their clones, contributing $8b_0$ to $b_1$. However, unlike the directed case, a transitive 3-cycle in $G_1$ can contain both a node in $G_0$ and its clone. Consider an arc $(u,v)$ of $G_0$. The oriented graph induced by the nodes $u,v$ and $u'$ or by the nodes $u,v$ and $v'$ forms a $T_3$ in $G_1$, as does the oriented graph induced by the nodes $u', v'$ and $u$, or by $u',v'$ and $v$. Thus, $b_1 = 8b_0+4e_0$.

Now, suppose that $s(1)=0$. Consider a transitive tournament on four nodes. Replacing none, all, or exactly one of the nodes with its clone results in a copy of $T_4$ in $G_1$. Replacing the source and its out-neighbor of largest out-degree, the sink and its in-neighbor of largest in-degree, or the two nodes other than the source and sink with their clones yields a copy of $T_4$. Finally, replacing all nodes except the source or sink node with their clones yields a copy of $T_4$, contributing a total of $11L_0$ to $L_1$. Replacing both the source and sink nodes, as well as possibly one other node, with their clones yields a Mixed tournament, contributing $3L_0$ to $m_1$. The remaining two cases result in one Loser tournament and one Winner tournament, contributing $L_0$ to each of $w_1$ and $\ell_1$.

Next, consider a Winner tournament in $G_0$. Replacing every node with its clone yields a Loser tournament, thereby contributing $w_0$ to $\ell_1$. Replacing none, exactly one, or exactly three of the nodes with their clones results in a Winner tournament, contributing $9w_0$ to $w_1$. Replacing the source node and any one other node with their clones yields a Mixed tournament, while replacing two non-source nodes with their clones yields a transitive tournament. Thus, $3w_0$ is contributed to $L_1$ and to $m_1$.

Next, consider a Loser tournament in $G_0$. As above, replacing every node with its clone yields a Winner tournament, contributing $\ell_0$ to $w_1$. Replacing none, exactly one, or exactly three of the nodes with their clones results in a Loser tournament, contributing $9\ell_0$ to $\ell_1$. Finally, replacing the sink node and any one other node with their clones yields a Mixed tournament, while replacing any two non-sink nodes with their clones yields a transitive tournament. Thus, $3\ell_0$ is contributed to $L_1$ and $m_1$.

Finally, consider a Mixed tournament in $G_0$. Replacing the two nodes with out-degree 2 by their clones yields a Winner tournament, contributing $m_0$ to $w_1$. Similarly, replacing the two nodes with out-degree one by their clones yields a Loser tournament, contributing $m_0$ to $\ell_1$. One of the two nodes of out-degree two has an in-neighbor with out-degree one. Replacing that pair of nodes with their clones yields a transitive tournament, as does replacing the two nodes of out-degree 2 and their common out-neighbor with their clones. Finally, replacing the two nodes of out-degree one and their common in-neighbor with their clones results in a transitive tournament. The remaining cases result in a Mixed tournament, contributing $3m_0$ to $L_1$ and $11m_0$ to $m_1$.

Let $a,b,c$ induce a copy of $D_3$ in $G_0$, with arcs $(a,b)$, $(b,c)$, and $(c,a)$. Among the four-node subsets containing at least one node from each of $\{a,a'\}$, $\{b,b'\}$, and $\{c,c'\}$, the sets $\{a,b,a',c'\}$, $\{b,c,a',b'\}$, and $\{a,c,b',c'\}$ induce Winner tournaments in $G_1$. The remaining nine subsets induce Mixed tournaments. This contributes $3a_0$ to $w_1$ and $9a_0$ to $m_1$.

For each copy of $D_3$ in $G_0$, replacing none, exactly one, or all three of the nodes with their clones creates a copy of $D_3$ in $G_1$. If exactly two nodes are replaced with their clones, then a copy of $T_3$ is obtained in $G_1$. This contributes $5a_0$ to $a_1$ and $3a_0$ to $b_1$.

Analogously, for each copy of $T_3$ in $G_0$, replacing only the source and sink nodes with their clones results in a copy of $D_3$ in $G_1$, while replacing any other subset of the three nodes with their clones results in a copy of $T_3$ in $G_1$. This contributes $b_0$ to $a_1$ and $7b_0$ to $b_1$.

Let $a,b,c$ be the nodes of a copy of $T_3$ in $G_0$, with arcs $(a,b)$, $(a,c)$, and $(b,c)$. Each of the sets $\{a,b,c,x\}$ with $x\in\{a',b',c'\}$, $\{x,b',c,c'\}$ with $x\in\{a,a'\}$, and $\{a',b,b',c\}$ induces a copy of $T_4$ in $G_1$. The subgraph induced by $\{a,b,b',c'\}$ is a Winner tournament. The subgraph induced by $\{a,a',b',c\}$ or by $\{a',b,c,c'\}$ is a Loser tournament. The remaining three subgraphs are Mixed tournaments. This contributes $6b_0$ to $L_1$, $b_0$ to $w_1$, $2b_0$ to $\ell_1$ and $3b_0$ to $m_1$.

Finally, consider an arc $(u,v)$ in $G_0$. The oriented graph induced by $\{u',v',u\}$ is a $D_3$, while those induced by $\{u,v,u'\}$, by $\{u,v,v'\}$ or by $\{u',v',v\}$ are all copies of $T_3$ in $G_1$. This contributes $e_0$ to $a_1$ and $3e_0$ to $b_1$.

The subgraph induced by $\{u,u',v,v'\}$ is a Loser tournament, contributing $e_0$ to $\ell_1$.

Summing the terms above yields the first six rows of the matrix when $s(1)=0$. The final two rows of both matrices follow from Theorem~\ref{prop:n and e}.
\end{proof}

Note that in the case that $G_0$ is a tournament, we have that $a_t + b_t = \binom{2^tn_0}{3}$ and $L_t + w_t + \ell_t + m_t = \binom{2^tn_0}{4}$.

Given tournaments $G$ and $H$ with $|V(G)|\geq|V(H)|$, let $n(H,G)$ be the number of subsets of $V(G)$ that induce a tournament isomorphic to $H$, and let $\mathrm{Aut}(H)$ be the automorphism group of $H$. We define the density of $H$ in $G$ by
$$d^*_G(H) = \frac{|\mathrm{Aut}(H)|n(H,G)}{|V(H)|!\binom{|V(G)|}{|V(H)|}}.$$
A sequence of tournaments $(G_n : n\geq 1)$ with $|V(G_n)|\to\infty$ is said to be \emph{quasirandom} if for every tournament $H$, we have that $$\lim_{n\to\infty} d^*_{G_n}(H)=2^{-\binom{|V(H)|}{2}}.$$
Quasirandom tournaments were first introduced in \cite{fan}. The uniformly random tournament on $n$ nodes labeled $1,2,\dots,n$ is obtained by independently orienting each pair of distinct nodes in either direction with probability $1/2$. In a quasirandom sequence, the density of any fixed tournament $H$ approaches its expected value in a uniformly random tournament. While calculating the density of every subtournament may appear challenging, it may be unnecessary. A tournament $G$ is \emph{quasirandom-forcing} if every sequence of tournaments $(G_n:n\geq1)$ satisfying $|V(G_n)|\to\infty$ and
$$\lim_{n\to\infty}d^*_{G_n}(G)=2^{-\binom{|V(G)|}{2}}$$
is quasirandom. We use the following quasirandom-forcing result.

\begin{theorem}[\cite{raz}]\label{thm:quasirandom forcing}
The linear order on four nodes is quasirandom-forcing.
\end{theorem}

Hence, to show that a sequence of tournaments is quasirandom, it is sufficient to consider the density of the linear order on four nodes. We take that approach below, correcting the proof of \cite[Theorem~6]{ILMT}.

\begin{theorem}\label{thm:0stepquasi}
    If $G_0$ is a tournament and $s$ is a generating sequence with infinite support, then $(G_t : t\geq 0)$ forms a quasirandom sequence.
\end{theorem}
\begin{proof}
We first consider the sequence consisting entirely of 0-entries. Recall that $L_t$ denotes the number of four-element subsets in $G_t$ that induce a linear order. By Theorem~\ref{thm:quasirandom forcing}, it suffices to show that the density of the linear order on four nodes in $G_t$ approaches $2^{-6}$ as $t\to \infty$ or, equivalently, that the proportion $L_t/\binom{|V(G_t)|}{4}$ approaches $\frac{3}{8}$ as $t \to \infty$.

From Theorem \ref{thm: motifs}, we have that
$$\iltovec{t} = \begin{bmatrix}
    11 & 3 & 3 & 3 & 0 & 6 & 0 & 0 \\
    1 & 9 & 1 & 1 & 3 & 1 & 0 & 0\\
    1 & 1 & 9 & 1 & 0 & 2 & 1 & 0\\
    3 & 3 & 3 & 11 & 9 & 3 & 0 & 0\\
    0 & 0 & 0 & 0 & 5 & 1 & 1 & 0\\
    0 & 0 & 0 & 0 & 3 & 7 & 3 & 0\\
    0 & 0 & 0 & 0 & 0 & 0 & 4 & 1\\
    0 & 0 & 0 & 0 & 0 & 0 & 0 & 2
\end{bmatrix}^t \iltovec{0}.$$
The above matrix $M$ can be diagonalized, so that $M = PDP^{-1}$, where
$$P = \begin{bmatrix}
    \phantom{-}3 & \phantom{-}72 & \phantom{-}33 & \phantom{-}3 & \phantom{-}3 & \phantom{-}3 & \phantom{-}6 & -3\\
    \phantom{-}1 & -24 & -25 & -3 & \phantom{-}0 & \phantom{-}0 & \phantom{-}0 & -1\\
    -1 & \phantom{-}24 & \phantom{-}23 & \phantom{-}0 & -3 & \phantom{-}0 & \phantom{-}0 & -1\\
    \phantom{-}3 & -72 & -75 & \phantom{-}0 & \phantom{-}0 & -3 & \phantom{-}0 & -3\\
    -2 & \phantom{-}48 & \phantom{-}48 & \phantom{-}0 & \phantom{-}0 & \phantom{-}0 & -1 & \phantom{-}0 \\
    -6 & -48 & \phantom{-}0 & \phantom{-}0 & \phantom{-}0 & \phantom{-}0 & -3 & \phantom{-}0 \\
    \phantom{-}12 & \phantom{-}0 & -48 & \phantom{-}0 & \phantom{-}0 & \phantom{-}0 & \phantom{-}0 & \phantom{-}0\\
    -24 & \phantom{-}0 & \phantom{-}0 & \phantom{-}0 & \phantom{-}0 & \phantom{-}0 & \phantom{-}0 & \phantom{-}0
\end{bmatrix} \text{ and } D=
\begin{bmatrix}
    2 & 0 & 0 & 0 & 0 & 0 & 0 & 0 \\
    0 & 4 & 0 & 0 & 0 & 0 & 0 & 0 \\
    0 & 0 & 4 & 0 & 0 & 0 & 0 & 0 \\
    0 & 0 & 0 & 8 & 0 & 0 & 0 & 0 \\
    0 & 0 & 0 & 0 & 8 & 0 & 0 & 0 \\
    0 & 0 & 0 & 0 & 0 & 8 & 0 & 0 \\
    0 & 0 & 0 & 0 & 0 & 0 & 8 & 0 \\
    0 & 0 & 0 & 0 & 0 & 0 & 0 & 16 \\
\end{bmatrix}.$$
By calculating the first entry of the vector
$$\iltovec{t} = PD^tP^{-1}\iltovec{0}$$ and dividing by $\binom{2^tn_0}{4}$, we have that
$$\lim_{t\to\infty} \frac{L_t}{\binom{2^tn_0}{4}} = \lim_{t\to\infty} \frac{1}{48\cdot 2^{4t}n_0^4} \bigg(27\cdot 2^{4t+4}\binom{n_0}{4} + 81\cdot 2^{4t+3}\binom{n_0}{3} + 63\cdot 2^{4t+2}\binom{n_0}{2} + 9\cdot 2^{4t+1}n_0\bigg) = \frac{3}{8}.$$

For a general generating sequence $s$ with infinite support, let $\delta_t=L_t/n_t^4-1/64$. The two matrices in Theorem~\ref{thm: motifs}, together with $L_t+w_t+\ell_t+m_t=\binom{n_t}{4}$, give
$$\delta_{t+1}=\begin{cases}
\frac12\delta_t+O(2^{-t}), & \text{if }s(t+1)=0,\\
\delta_t+O(2^{-t}), & \text{if }s(t+1)=1.
\end{cases}$$
For any fixed $N$, infinitely many subsequent 0-steps force the contribution from $\delta_N$ to tend to zero, while the sum of the absolute values of the later error terms is at most $C\sum_{j=N}^{\infty}2^{-j}$ for a constant $C$ independent of $N$. Letting $N\to\infty$ gives $\delta_t\to0$, and hence, $L_t/\binom{n_t}{4}\to3/8$, as required.
\end{proof}

We conclude the section by noting that the ILTO process on a tournament yields, in a certain sense, continuum-many pairwise non-isomorphic sequences of tournaments.  Hence, combined with Theorem~\ref{thm:0stepquasi}, the ILTO model yields continuum-many non-isomorphic sequences of quasirandom tournaments.

\begin{theorem}
Given a tournament $G_0$ and two distinct generating sequences $s$ and $s'$ with $s(1)=s'(1)$, there is some $t\geq1$ such that $\ILTO_{t,s}(G_0)$ is not isomorphic to $\ILTO_{t,s'}(G_0)$.
\end{theorem}

\begin{proof} Let $G_t=\ILTO_{t,s}(G_0)$ and $G_t'=\ILTO_{t,s'}(G_0)$. Let $t_0+1$ be the first index at which $s$ and $s'$ differ. Since $s(1)=s'(1)$, we have $t_0>0$, and the two constructions agree through time $t_0$. Interchanging $s$ and $s'$ if necessary, assume that $s(t_0+1)=1$ and $s'(t_0+1)=0$.

Let $r$ be the number of nodes in $G_{t_0}$  with out-degree $2^{t_0-1}n_0$. A node with out-degree $\delta$ at time-step $t_0$ must have double the out-degree, $2\delta$, at time-step $t_0 + 1$. Every clone node in $G'_{t_0 + 1}$ has $2^{t_0}n_0 - 1$ out-neighbors other than its parent. Thus, $G_{t_0+1}'$ has $r+2^{t_0}n_0$ nodes with out-degree $2^{t_0}n_0$. However, note that $G_{t_0+1}$ has only $r$ nodes with out-degree $2^{t_0}n_0$, since all clones generated by a 1-step have odd out-degree. Altogether, we have that the tournaments $G_{t_0+1}$ and $G_{t_0+1}'$ are not isomorphic, since they do not have the same out-degree distribution. That is, the desired result holds for $t=t_0+1$.
\end{proof}

\section{Graph Parameters}

\subsection{Domination number}

For an oriented graph $G$, we say that $S\subseteq V(G)$
 is an \emph{in-dominating} set for $G$ if every node in $V(G)\setminus S$ has an in-neighbor in $S$. Out-dominating sets are defined analogously. An in-dominating (respectively, out-dominating) set for $G$ is \emph{minimal} if it does not contain any in-dominating (respectively, out-dominating) set as a proper subset. The \emph{in-domination number} of $G$, denoted by $\gamma^-(G)$, is the smallest integer $k$ such that $G$ has an in-dominating set of order $k$. The \emph{out-domination number} of $G$, denoted by $\gamma^+(G)$, is defined analogously. If $G$ is a tournament, then $\gamma^-(G) = 1$ if and only if $G$ has a source node. Similarly, for a tournament $G$, $\gamma^+(G) = 1$ if and only if $G$ has a sink node. 

 We first show that both the in-domination number and out-domination number of an ILTO graph are constant across 1-steps. 
 
\begin{theorem}
\label{prop:dom 1step}
    Let $G_0$ be an oriented graph and $G_1$ the oriented graph obtained via a 1-step. We then have that $\gamma^-(G_1) = \gamma^-(G_0)$ and $\gamma^+(G_1) = \gamma^+(G_0)$. 
\end{theorem}
\begin{proof}
    Let $S$ be an in-dominating set for $G_0$ and let $S' = \{x' : x\in S\}$. In $G_1$, each $s\in S$ is in-dominated by $s'\in S'$. If $x\in V(G_0)\setminus S$, then $x$ is in-dominated in $G_0$ by some $s\in S$. In $G_1$, both $x$ and its clone $x'$ are in-dominated by $s'\in S'$. Thus, $S'$ is an in-dominating set for $G_1$ and $\gamma^-(G_1) \leq \gamma^-(G_0)$. 

    Conversely, let $S$ be an in-dominating set in $G_1$ and let $S^*$ be the set obtained from $S$ by replacing any $v\in V(G_1)\setminus V(G_0)$ by its parent node in $G_0$. For $y\in V(G_0)\setminus S^*$, if $y$ is in-dominated in $G_1$ by $s\in S\cap V(G_0)$, then $y$ is also in-dominated by $s$ in $G_0$. Otherwise, $y$ is in-dominated in $G_1$ by some $s\in S$ with $s\not\in V(G_0)$. The parent node $s^*$ of $s$ in $G_0$ is in the set $S^*$ and in-dominates $y$ in $G_0$. Thus, $S^*$ is in-dominating in $G_0$ and $\gamma^-(G_0)\leq \gamma^-(G_1)$. 

    Similarly, let $S$ be an out-dominating set for $G_0$. In $G_1$, for each $x\in S$, its clone node $x'$ is out-dominated by $x$. For each $y\in V(G_0)\setminus S$, there exists $x\in S$ such that $(y,x)$ is an arc. We then have that $(y',x)$ is an arc in $G_1$, so both $y$ and $y'$ are out-dominated by a node of $S$. Thus, $\gamma^+(G_1)\leq \gamma^+(G_0)$. 

    Similarly, let $S$ be an out-dominating set for $G_1$. Let $S^*$ be the set obtained from $S$ by replacing any node $x'\in S\cap(V(G_1)\setminus V(G_0))$ by its parent node in $G_0$, so that $S^*\subseteq V(G_0)$. Each node $y\in V(G_0)\setminus S^*$ is out-dominated in $G_1$ by some $s\in S$. If $s\in V(G_0)$, then $y$ is also out-dominated by $s$ in $G_0$. Otherwise, $s$ was replaced by its parent node, $s^*$, in $S^*$. It follows that $y$ is out-dominated by $s^*$ in $G_0$. Thus, $\gamma^+(G_1)=\gamma^+(G_0)$. 
\end{proof}

That the out-domination number is unchanged by a 0-step follows from the same argument as in the proof of Theorem~\ref{prop:dom 1step}. Thus, we omit the proof.  

\begin{lemma}
    Let $G_0$ be an oriented graph and $G_1$ the oriented graph obtained via a 0-step. We then have that $\gamma^+(G_1) = \gamma^+(G_0)$. 
\end{lemma}

If $G$ is an oriented graph with an in-dominating set (respectively, out-dominating set) $S\subseteq V(G)$ and $H$ is the reverse of $G$, then $S$ is an out-dominating (respectively, in-dominating) set for $H$. Recall that if $G_1$ is constructed from $G_0$ via a 0-step, then the set of clones in $G_1$ induces a copy of the reverse of $G_0$. We can use these facts to obtain an in-dominating set for $G_1$ containing only clone nodes. 

\begin{theorem}
\label{thm:in and out dominating}
    Let $G_0$ be an oriented graph and $G_1$ the oriented graph obtained via a 0-step. For $S\subseteq V(G_0)$, let $S' = \{x' : x\in S\}$. The set $S'$ is an in-dominating set for $G_1$ if and only if $S$ is both an in- and out-dominating set for $G_0$. 
\end{theorem}

\begin{proof}
    Suppose that $S'$ is an in-dominating set for $G_1$. For every $x\in V(G_0)\setminus S$, there exist nodes $u',v'\in S'$ such that $(u',x)$ and $(v',x')$  are arcs in $G_1$. By the construction of a 0-step, $(u,x)$ and $(x,v)$ are arcs in $G_0$. Thus, $x$ is both in-dominated and out-dominated by a node in $S$. 

    Conversely, if $S$ is both in-dominating and out-dominating for $G_0$, then for every $x\in V(G_0)\setminus S$, there exist nodes $u,v\in S$ such that $(u,x)$ and $(x,v)$ are arcs in $G_0$. In $G_1$, we have arcs $(u',x)$ and $(v',x')$, so both $x$ and $x'$ are in-dominated by a node in $S'$. Moreover, every $x\in S$ is in-dominated by the node $x'\in S'$. Thus, $S'$ is an in-dominating set for $G_1$.  
\end{proof}

In \cite[Lemma~1 (2)]{ILMT}, it was claimed that if $G_1$ is a tournament obtained from the tournament $G_0$ via a 0-step and $S$ is a minimal in-dominating set for $G_0$, then $S'=\{x':x\in S\}$ is an in-dominating set for $G_1$. Moreover, \cite[Theorem~8]{ILMT} asserted that $S'$ is a minimal in-dominating set for $G_1$, and \cite[Corollary~3]{ILMT} claimed that $\gamma^-(G_1)=\gamma^-(G_0)$ when $G_0$ has no source. The following example disproves all three claims.

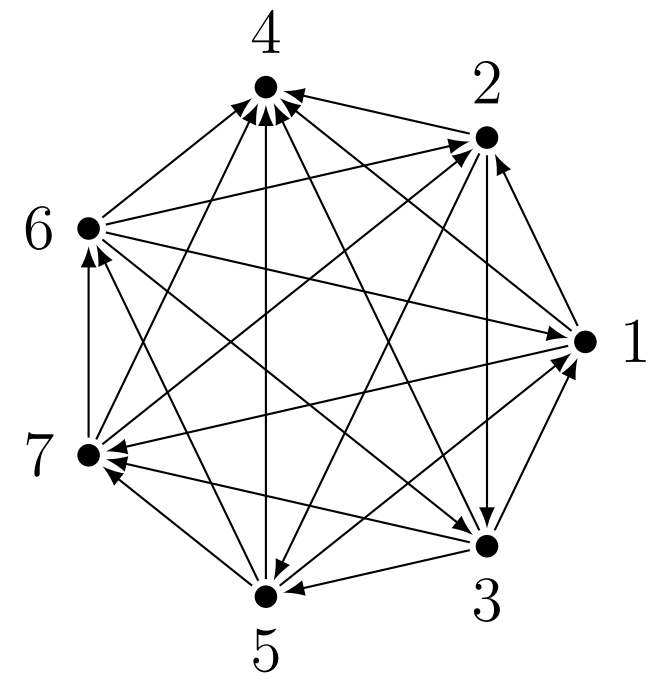
\begin{figure}[ht]
    \centering
   \begin{tikzpicture}[
    vertex/.style={circle, fill=black, inner sep=1.5pt},
    edge/.style={->, >={latex[length=2.5mm, width=1.5mm]}, shorten >=1pt, shorten <=1pt}
    ]

    \node[vertex, label=right:1] (1) at (1.5*1.1524, 0) {};
    \node[vertex, label=above:2] (2) at (1.5*0.7185, 1.5*0.901) {};
    \node[vertex, label=below:3] (3) at (1.5*0.7185, 1.5*-0.901) {};
    \node[vertex, label=above:4] (4) at (1.5*-0.2564, 1.5*1.1235) {};
    \node[vertex, label=below:5] (5) at (1.5*-0.2564, 1.5*-1.1235) {};
    \node[vertex, label=left:6] (6) at (1.5*-1.0383, 1.5*0.5) {};
    \node[vertex, label=left:7] (7) at (1.5*-1.0383, 1.5*-0.5) {};

    \path[edge] (1) edge (2);
    \path[edge] (1) edge (4);
    \path[edge] (1) edge (7);
    \path[edge] (2) edge (3);
    \path[edge] (2) edge (4);
    \path[edge] (2) edge (5);
    \path[edge] (3) edge (1);
    \path[edge] (3) edge (4);
    \path[edge] (3) edge (5);
    \path[edge] (3) edge (7);
    \path[edge] (5) edge (1);
    \path[edge] (5) edge (4);
    \path[edge] (5) edge (6);
    \path[edge] (5) edge (7);
    \path[edge] (6) edge (1);
    \path[edge] (6) edge (2);
    \path[edge] (6) edge (3);
    \path[edge] (6) edge (4);
    \path[edge] (7) edge (2);
    \path[edge] (7) edge (4);
    \path[edge] (7) edge (6);

    \end{tikzpicture}
    \caption{A tournament $G_0$ with $\gamma^-(G_0)=2$ and $\gamma^-(G_1)=3$, where $G_1$ is constructed via a 0-step.}
    \label{fig:counterex domination}
    \end{figure}

    Let $G_0$ be the tournament in Figure \ref{fig:counterex domination}. Since $G_0$ has no source node and $S = \{5,6\}$ is an in-dominating set for $G_0$, we have that $\gamma^-(G_0)=2$. However, in $G_1$, the set $S'=\{5',6'\}$ does not in-dominate either of the clones $1'$ or $4'$. Note that both of the corresponding parent nodes, $1$ and $4$, are out-neighbors of every node in $S$. It is straightforward to check that there is no in-dominating set of order two for $G_1$, but $\{4',5,6\}$ in-dominates every node of $G_1$, so $\gamma^-(G_1)=3$. We also note that $\{4\}$ is an out-dominating set, so $S_1 = \{4,5,6\}$ is both in-dominating and out-dominating in $G_0$ and, by Theorem~\ref{thm:in and out dominating}, $S_1' = \{4',5',6'\}$ is in-dominating in $G_1$. 

    We can always take the union of an in-dominating set and an out-dominating set in $G_0$ and clone every node to obtain an in-dominating set for $G_1$. Thus, we have the following upper bound on the in-domination number for an ILTO graph constructed via a 0-step. 

\begin{corollary}
\label{cor:in plus out}
    Let $G_0$ be an oriented graph and $G_1$ the oriented graph obtained via a 0-step. We have that $$\gamma^-(G_1)\leq \gamma^-(G_0)+\gamma^+(G_0).$$
\end{corollary}

To obtain an upper bound on the in-domination number for $G_t = \ILTO_{t,s}(G_0)$, we apply Theorem~\ref{prop:dom 1step} for each 1-entry in $s$ and Corollary \ref{cor:in plus out} for each 0-entry. Thus, the in-domination number increases by at most $\gamma^+(G_0)$ for each 0-entry in $s$. 

\begin{corollary}
    Let $G_0$ be an oriented graph and $G_t = \ILTO_{t,s}(G_0)$. For $t\geq 1$, we have that $$\gamma^-(G_t)\leq \gamma^-(G_0) + k\gamma^+(G_0),$$
    where $k = |\{1\leq i \leq t : s(i)=0\}|$. 
\end{corollary}

The following bound is independent of the generating sequence and the number of time-steps.

\begin{theorem}\label{thm:bounded indomination}
Let $G_0$ be a connected oriented graph, let $s$ be a generating sequence, and let $G_t=\ILTO_{t,s}(G_0)$. For every $t\geq 0$, we have that
$$\gamma^-(G_t)\leq 2\bigl(\gamma^-(G_0)+\gamma^+(G_0)\bigr).$$
\end{theorem}

\begin{proof}
If $s$ contains no 0-entry, then Theorem~\ref{prop:dom 1step} gives $\gamma^-(G_t)=\gamma^-(G_0)$ for all $t$, so the bound holds. Suppose first that $G_0$ has at least two nodes, and let $q$ be the index of the first 0-entry. By Theorem~\ref{prop:dom 1step} and Corollary~\ref{cor:in plus out}, we have that
$$\gamma^-(G_q)\leq\gamma^-(G_0)+\gamma^+(G_0).$$
The bound also holds for $t<q$, since all preceding steps are 1-steps.

We claim that $G_q$ has no source. Every parent has its clone as an in-neighbor. For each clone $v'$, choose a neighbor $u$ of $v$ in $G_{q-1}^*$, which exists since $G_{q-1}$ is connected and has at least two nodes. If $u\to v$, then $u\to v'$; if $v\to u$, then $u'\to v'$. Thus, every clone also has an in-neighbor.

Let $S$ be a minimum in-dominating set of $G_q$. For each node of $S$ with no in-neighbor in $S$, choose an in-neighbor in $G_q$ and add it to $S$, obtaining a set $T$. Every added node lies outside $S$ and therefore has an in-neighbor in $S$. Consequently, every node of $G_q$, including every node of $T$, has an in-neighbor in $T$, and $|T|\leq 2|S|$.

This property persists through either type of time-step while keeping $T$ fixed. Indeed, all existing arcs remain, and if $u\in T$ is an in-neighbor of a node $v$, then $u$ is also an in-neighbor of its new clone $v'$. Hence, $T$ is in-dominating in every $G_t$ with $t\geq q$, giving
$$\gamma^-(G_t)\leq |T|\leq 2\gamma^-(G_q)\leq 2\bigl(\gamma^-(G_0)+\gamma^+(G_0)\bigr).$$

Finally, if $G_0$ consists of a single node, then $G_1$ consists of a single arc and $\gamma^-(G_1)=\gamma^+(G_1)=1=\gamma^-(G_0)=\gamma^+(G_0)$. Applying the case just proved with $G_1$ as the base graph gives the bound for all $t\geq 1$, and the case $t=0$ is immediate.
\end{proof}

As demonstrated by the above example, if $S\subseteq V(G_0)$ is an in-dominating set for $G_0$ and $G_1$ is constructed via a 0-step, then $S' = \{x' : x\in S\}$ is not necessarily an in-dominating set for $G_1$. The nodes outside $S'$ that have no in-neighbor in $S'$ are exactly the clones of nodes in $V(G_0)\setminus S$ that have no out-neighbor in $S$, since a 0-step reverses arcs between clones. We use this information to characterize the in-dominating sets of $G_0$ whose clones form in-dominating sets in $G_1$. 

\begin{lemma}
\label{lem:indom}
    Let $G_0$ be an oriented graph and $G_1$ the oriented graph obtained via a 0-step. Let $S$ be an in-dominating set for $G_0$. We then have that $S' = \{x' : x\in S\}$ is an in-dominating set for $G_1$ if and only if $N^+_{G_0}(v)\cap S \neq \emptyset$ for all $v\in V(G_0)\setminus S$. 
\end{lemma}

\begin{proof}
    Suppose that there exists $v\in V(G_0)\setminus S$ such that $N^+_{G_0}(v)\cap S = \emptyset$ and consider the arcs between the node $v'$ and the set $S'$. For all $s\in S$, we either have the arc $(s,v)$ or no arc between $s$ and $v$. If the arc $(s,v)$ exists in $G_0$, then we have the arc $(v',s')$ in $G_1$. If there is no arc between $s$ and $v$ in $G_0$, then there is no arc between $s'$ and $v'$ in $G_1$. Thus, $v'$ has no in-neighbor in $S'$ and so $S'$ is not an in-dominating set. 

    Suppose instead that $N^+_{G_0}(v)\cap S \neq \emptyset$ for all $v\in V(G_0)\setminus S$. If $x\in S$, then $x'\in S'$ is an in-neighbor of $x$. If $x \in V(G_0)\setminus S$, then $x$ has some in-neighbor in $S$, say $s$. Hence, $s'\in S'$ is an in-neighbor of $x$. Finally, consider $x'$, where $x\in V(G_0)\setminus S$. There exists $s\in N^+_{G_0}(x)\cap S$. Since $(x,s)$ is an arc in $G_0$, we must have the arc $(s',x')$ in $G_1$. Thus, $S'$ is an in-dominating set in $G_1$. 
\end{proof}

Suppose that $\{v\}$ is an in-dominating set for $G_0$, where $G_0$ has at least two nodes, and let $G_1$ be obtained from
$G_0$ via a 0-step. Since $v$ has an arc to every other node of $G_0$, the clone $v'$ has no arc to any other clone.
Thus, $\{v'\}$ is not in-dominating in $G_1$, but $\{v,v'\}$ is. Every parent has its clone as an in-neighbor,
and every clone other than $v'$ has $v$ as an in-neighbor. Further, $v'$ has $u'$ as an in-neighbor for any
$u\in V(G_0)\setminus\{v\}$. Hence, no singleton is in-dominating in $G_1$, and we obtain
the following corollary.

\begin{corollary}\label{cor:no source}
Let $G_0$ be an oriented graph of order at least two and $G_1$ the oriented graph obtained via a 0-step. If $\gamma^-(G_0)=1$, then $\gamma^-(G_1)=2$.
\end{corollary}

When $G_0$ is a tournament and $S$ is an in-dominating set for $G_0$, every node $v\in V(G_0)\setminus S$ is in either the in-neighborhood or the out-neighborhood of a node $s\in S$. To determine whether $S'$ is in-dominating for $G_1$, we must determine whether each node $v \in V(G_0)\setminus S$ is in the in-neighborhood of some $s\in S$ or, equivalently, whether the intersection of the out-neighborhoods of the nodes in $S$ is empty. 

\begin{corollary}
    Let $G_0$ be a tournament and $G_1$ the tournament obtained via a 0-step. Let $S$ be an in-dominating set for $G_0$. We then have that $S' = \{x' : x\in S\}$ is an in-dominating set for $G_1$ if and only if $\bigcap_{s\in S} N^+_{G_0}(s) = \emptyset$. 
\end{corollary}

\begin{proof}
    If $S'$ is an in-dominating set for $G_1$, then by Lemma~\ref{lem:indom}, $N^+_{G_0}(v)\cap S \neq \emptyset$ for all $v\in V(G_0)\setminus S$. This means that for each $v\in V(G_0)\setminus S$, there is some $s\in S$ such that $s\in N^+_{G_0}(v)$. Thus, $v\not\in N^+_{G_0}(s)$ and so $v\not\in \bigcap_{s\in S} N^+_{G_0}(s)$. Since $\bigcap_{s\in S} N^+_{G_0}(s)\subseteq V(G_0)\setminus S$, we have that $\bigcap_{s\in S} N^+_{G_0}(s)=\emptyset$. 

    Suppose now that $\bigcap_{s\in S} N^+_{G_0}(s)=\emptyset$ and consider $v\in V(G_0)\setminus S$. Since $v\not\in \bigcap_{s\in S} N^+_{G_0}(s)$, there exists $s\in S$ such that $v\not\in N^+_{G_0}(s)$. Since $G_0$ is a tournament, there is an arc between $v$ and $s$, which must go from $v$ to $s$. Therefore, $s\in N^+_{G_0}(v)\cap S$. It then follows from Lemma~\ref{lem:indom} that $S'$ is an in-dominating set for $G_1$. 
\end{proof}

We make use of the following result of Erd\H{o}s \cite{erdos} to improve the upper bound on the in-domination number of an ILTO graph when the base graph is a tournament. 

\begin{lemma}
\label{lem:erdos}
    \cite{erdos} For every tournament $G$ of order $n$, $\gamma^-(G)\leq \lceil \log_2(n+1)\rceil$. 
\end{lemma}

Now, suppose that $S$ is an in-dominating set for the tournament $G_0$ and $G_1$ is obtained from $G_0$ via a 0-step. If $S'$ is not an in-dominating set for $G_1$, then the nodes in $V(G_1)\setminus S'$ with no in-neighbor in $S'$ induce a nonempty subtournament $H$ of $G_1$. Instead of adding the clones of the nodes of an out-dominating set of $G_0$ to $S'$ to obtain an in-dominating set for $G_1$, we can add the nodes of an in-dominating set for the tournament $H$. By the above lemma, there exists such an in-dominating set containing at most $\lceil \log_2(|V(H)|+1)\rceil$ nodes. 

\begin{theorem}
    Let $G_0$ be a tournament and $G_1$ the tournament obtained via a 0-step. For an in-dominating set $S\subseteq V(G_0)$, let $N_S = \bigcap_{s\in S} N^+_{G_0}(s)$. 
    We then have that
    $$\gamma^-(G_1)\leq \gamma^-(G_0) + \min\{\gamma^+(G_0),n^*\},$$
    where $n^* = \min_{S\in \mathcal{S}} \lceil \log_2\left(|N_S|+1\right)\rceil$ and $\mathcal{S}$ is the set of all in-dominating sets of $G_0$ of cardinality $\gamma^-(G_0)$. 
\end{theorem}
\begin{proof}
    Choose $S\in\mathcal{S}$ such that $\lceil\log_2(|N_S|+1)\rceil=n^*$. Thus, $S'$ in-dominates every node outside $S'\cup N_S'$. If $N_S=\emptyset$, then $S'$ is in-dominating in $G_1$ and $n^*=0$, so the claimed bound follows. We may therefore assume that $N_S\neq\emptyset$. The nodes in $N_S'$ induce a tournament of order $|N_S|$ which, by Lemma~\ref{lem:erdos}, has an in-dominating set $S_{N_S}$ of cardinality at most $\lceil \log_2 \left(|N_S|+1\right)\rceil$. Hence, the set $S'\cup S_{N_S}$ is an in-dominating set for $G_1$ of order at most $\gamma^-(G_0)+\lceil \log_2\left(|N_S|+1\right)\rceil$. The result follows from Corollary \ref{cor:in plus out}. 
\end{proof}

The tournament in Figure \ref{fig:counterex domination} shows that this upper bound is tight. 

\subsection{Cop number}

\emph{Cops and Robber} is a game that can be played on an oriented graph $G$. One player controls a set of cops, and the other controls a single robber. The cops first choose their initial nodes, after which the robber chooses an initial node. Multiple cops may occupy the same node. Both players know the positions of all cops and the robber throughout the game. Play then proceeds in rounds, with the cops moving first and the robber moving second. Each cop and the robber may move from their current node $u$ along an arc $(u,v)$ to $v$. Players can pass or remain on their own node. Any subset of cops may move in a given round. The cops win the game if, after a finite number of rounds, one of them occupies the same node as the robber. This situation is called a \emph{capture}. The robber wins if they can evade capture indefinitely. 

The minimum number of cops required to win is a well-defined positive integer, called the \emph{cop number} of $G$, denoted by $c(G)$. For additional background on the cop number of a graph, see the book \cite{cops}. 

The cop number in the ILTO model does not change with a 1-step, as the following argument shows. First, we show that the cop number does not increase following a 1-step. 

\begin{lemma}
If $G_0$ is an oriented graph and $G_1$ the oriented graph obtained via a 1-step, then $c(G_1) \leq c(G_0).$
\end{lemma}
\begin{proof}
Let $c(G_0)$ cops start by exclusively moving among the clone nodes. The clone subgraph is a copy of $G_0$, and mapping each parent to its clone turns every robber move into a legal move or a pass in this copy. They first capture the robber's image in the subgraph induced by the clone nodes, meaning that a cop occupies $u'$ while the robber occupies either $u$ or $u'$. 
We can therefore assume that at some stage, the robber is on some node $v$ and a cop is on $v'$. 
If the robber moves to a node $w$ along some arc $(v,w)$, then there exists the arc $(v',w)$.  The cop on $v'$ then moves to $w$ and captures the robber. If the robber instead moves from $v$ to a clone $w'$, then $(v,w)$ is an arc in $G_0$, so $(v',w')$ is an arc in $G_1$, and the cop captures the robber there.
If the robber remains on the node $v$, then the cop on $v'$ can use the arc $(v',v)$ to move to $v$ and capture the robber.  
\end{proof}

For a graph homomorphism $\phi:V(G_1)\to V(G_2)$ and $S \subseteq V(G_1)$, we define $\phi(S) = \{\phi(s) : s \in S\}$. 
A homomorphism $\phi:V(G_1)\to V(G_2)$ is  \emph{full} if for each pair of nodes $\{u,v\}\subseteq V(G_1)$ where neither the arc $(u,v)$ nor $(v,u)$ occurs in $G_1$, the arcs $(\phi(u),\phi(v))$ and $(\phi(v),\phi(u))$ do not occur in $G_2$.

\begin{theorem}[\cite{ILMT}] \label{thm:full_homomorphisms} 
If $G,H$ are oriented graphs and there is a surjective full homomorphism $\phi:V(G)\to V(H)$, then $c(G)\geq c(H)$.
\end{theorem}
\begin{proof}
The robber can escape capture on $H$ when $c(H)-1$ cops are playing, and so let $\mathcal{R}_H(C,R)=R'$ be such a robber strategy, where $\mathcal{R}_H(C,\emptyset)$ indicates where the robber first appears. We use this robber strategy to avoid capture in $G$ against any strategy $\mathcal{C}_G$ for $c(H)-1$ cops. Suppose the cops initially occupy the set $C_0$ in $G$. The robber chooses an initial node $R_0\in\phi^{-1}(\mathcal{R}_H(\phi(C_0),\emptyset))$, which exists by surjectivity.

In the first round, the cops move from $C_0$ to $C_1=\mathcal{C}_G(C_0,R_0)$. Since $\phi$ preserves arcs, their images make legal moves in $H$, and the robber's strategy ensures that $\phi(R_0)\notin\phi(C_1)$. Thus, the cops do not capture the robber at $R_0$. If the robber's strategy in $H$ leaves it at $\phi(R_0)$, then set $R_1=R_0$, and the robber remains there. Otherwise, choose $R_1\in\phi^{-1}(\mathcal{R}_H(\phi(C_1),\phi(R_0)))$, which exists by surjectivity. Fullness implies that $R_0$ and $R_1$ are adjacent. Since $\phi$ preserves arcs and $H$ is oriented, this arc must be $(R_0,R_1)$. The robber moves to $R_1$. Since $\phi(R_1)\notin\phi(C_1)$, we have $R_1\notin C_1$.

We repeat this process in each subsequent round, with the robber avoiding capture. Thus, $c(H)-1$ cops are insufficient to capture the robber in $G$, and the result follows.
\end{proof}

We demonstrate next that the cop number of an ILTO graph also cannot decrease following a 1-step. 

\begin{lemma}
If $G_0$ is an oriented graph and $G_1$ is the oriented graph obtained via a 1-step, then
$c(G_1) \geq c(G_0).$
\end{lemma}
\begin{proof}
Let $\phi:V(G_1)\to V(G_0)$ be the function that maps clones to their parents and leaves parents unchanged; that is, $\phi(u')=u$ and $\phi(u)=u$. Each move in $G_1$ maps under $\phi$ to a legal move or a pass in $G_0$. Thus, a winning strategy for $c(G_1)$ cops in $G_1$, with the robber restricted to parent nodes, gives a winning strategy in $G_0$ by mapping each cop's position under $\phi$. Hence, $c(G_0)\leq c(G_1)$.
\end{proof}

\begin{corollary}
    Let $G_0$ be an oriented graph and $G_1$ the oriented graph obtained via a 1-step. We then have that $c(G_1) = c(G_0)$. 
\end{corollary}

The in-domination number provides an upper bound for the cop number of an oriented graph. Together with Theorem~\ref{thm:bounded indomination}, this gives an upper bound on the cop number depending only on the base graph, for every generating sequence.

For an oriented graph $G$, let $\alpha(G)$ denote the independence number of the underlying graph of $G$. We show that the independence number of the base graph can be used to bound the cop number of the ILTO graph obtained via a 0-step. 

\begin{theorem}\label{thm: cops}
    Let $G_0$ be an oriented graph and $G_1$ the oriented graph obtained via a 0-step. We then have that $c(G_1)\leq \min\{3\alpha(G_0),2\alpha(G_0)+c(G_0)\}$. 
\end{theorem}

\begin{proof}
Let $S$ be an independent set of cardinality $\alpha(G_0)$ in $G_0^*$. For each $s\in S$, place one cop on the node $s$ and two cops on the clone node $s'$. The robber cannot occupy a node of $S\cup S'$, since each such node is occupied by a cop. We will show that the cops can capture the robber in at most two moves.

First, suppose that the robber starts on a clone node $v'$ in $G_1$. Since $S$ is a maximal independent set in $G_0^*$, there exists a node $s\in S$ such that there is an arc between $v$ and $s$ in $G_0$. If $(s,v)$ is an arc in $G_0$, then the cop on the node $s$ can move along the arc $(s,v')$ in $G_1$ to capture the robber. Otherwise, $(v,s)$ is an arc, and one of the cops on $s'$ can move along the arc $(s',v')$ to capture the robber.

Suppose instead that the robber starts on a parent node $v$. As above, there is some node $s\in S$ such that there is an arc between $v$ and $s$ in $G_0$. If $(s,v)$ is an arc, then the cop on $s$ can move along this arc to capture the robber. Otherwise, $(v,s)$ is an arc in $G_0$, so the cops on $s$ and $s'$ cannot capture the robber immediately. One of the two cops on $s'$ will move along the arc $(s',v')$ to the node $v'$. Note that there is still at least one cop occupying every node in $S$ and its clone. If the robber moves onto a node of $S\cup S'$, then capture is immediate. Thus, we may assume that any new position lies outside $S\cup S'$. Now, if the robber does not move, then the cop on $v'$ can move along the arc $(v',v)$ to capture the robber. If the robber moves along an arc $(v,w')$, then there exists a node $z\in S$ such that either $(z,w)$ or $(w,z)$ is an arc in $G_0$. In the former case, the cop on $z$ moves to capture the robber; in the latter case, one of the cops on $z'$ moves to capture the robber. Finally, if the robber moves along an arc $(v,w)$, then the cop on $v'$ can move along the arc $(v',w)$ to capture the robber. Thus, $c(G_1)\leq 3\alpha(G_0)$.

Suppose instead that the cops begin by placing one cop on each node $s$ and one cop on its clone $s'$, for each $s\in S$, and placing $c(G_0)$ cops in $G_0$ according to a winning cop strategy in $G_0$. The $2\alpha(G_0)$ cops on the nodes of $S\cup S'$ will pass unless one of them can move along an arc to immediately capture the robber. As above, if the robber starts on a clone node $v'$ in $G_1$, then the robber will be captured. Similarly, if the robber ever moves to occupy a clone node, then one of the $2\alpha(G_0)$ cops on $S\cup S'$ will capture the robber. If the robber remains on the nodes in $V(G_0)$, then the $c(G_0)$ cops have a winning strategy to capture the robber. Thus, $c(G_1)\leq 2\alpha(G_0)+c(G_0)$.
\end{proof}

Note that if $G_0$ is a tournament, then its underlying graph is a complete graph, which has independence number 1. 

\begin{corollary}
    [\cite{ILMT}] Let $G_0$ be a tournament and $G_1$ the tournament obtained via a 0-step. We then have that $c(G_1)\leq 3$. 
\end{corollary}

We note that it is unknown whether there are any oriented graphs that attain the bound in Theorem~\ref{thm: cops}, even when the base graph is restricted to tournaments. The following elementary lemma is used to show that the cop number after a 0-step is at least two whenever the base graph has at least two nodes.

\begin{lemma}[\cite{Halliday_thesis}] Let $G$ be an oriented graph. If $G$ has no source node, then $c(G)\geq 2$.   
\end{lemma}

If $G_0$ has at least two nodes, then every node of $G_0$ has a neighbor, since $G_0$ is connected. In $G_1$, every parent
has its clone as an in-neighbor. For each clone $v'$, a neighbor $u$ of $v$ gives an in-neighbor of $v'$: this is $u$
if $(u,v)$ is an arc in $G_0$, and $u'$ if $(v,u)$ is an arc. Thus, $G_1$ has no source.

\begin{theorem}
    If $G_0$ is an oriented graph of order at least two and $G_1$ is the oriented graph obtained via a 0-step, then $c(G_1) \geq 2$. 
\end{theorem}

\subsection{Dichromatic number}

We finish the section by considering dicolorings of ILTO oriented graphs. A $k$-\emph{dicoloring} of an oriented graph $G$ is a function $\varphi:V(G) \to \{1,2,\dots,k\}$ such that the subgraph induced by the set $\{v\in V(G) : \varphi(v)=i\}$ is an acyclic oriented graph, for $1\leq i \leq k$. The \emph{dichromatic number} of an oriented graph $G$, denoted by $\vec{\chi}(G)$, is the minimum $k$ such that $G$ admits a
$k$-dicoloring. The dichromatic number was introduced by Erd\H{o}s \cite{erdos2} and Neumann-Lara \cite{neumann} and rediscovered by Mohar \cite{mohar2} and by Bokal et al.\ \cite{mohar1}. 

For a 1-step, the dichromatic number remains unchanged.

\begin{theorem}
 If $G_0$ is an oriented graph and $G_1$ is the oriented graph obtained via a 1-step, then $\vec\chi(G_1) = \vec\chi(G_0)$. 
\end{theorem}

\begin{proof} Let $S=\{v_1,v_2,\dots,v_j\}$ be a color class in a $\vec\chi(G_0)$-dicoloring of $G_0$. Since $G_0[S]$ is acyclic, we may label its nodes so that every arc $(v_\ell,v_m)$ with both endpoints in $S$ satisfies $\ell<m$. Let $v'_1,v'_2,\dots,v'_j$ be the clones of $v_1,v_2,\dots,v_j$, respectively. It follows that $v'_1,v_1,v'_2,v_2,\dots,v'_j,v_j$ forms an acyclic oriented subgraph of $G_1$. 

Repeating this for each color class, we obtain a $\vec\chi(G_0)$-dicoloring of $G_1$. Hence, $\vec\chi(G_1) \leq \vec\chi(G_0)$. Since $G_0$ is a subgraph of $G_1$, we also have that $\vec\chi(G_1) \geq \vec\chi(G_0)$. 
\end{proof}

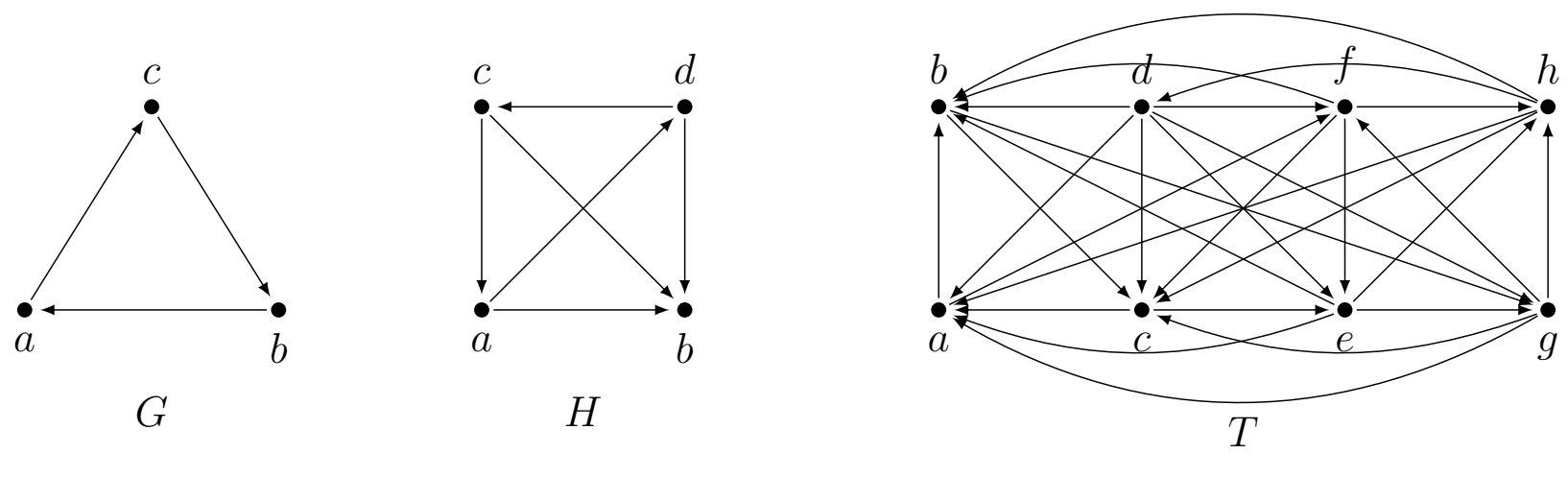
\begin{figure}[ht]
    \centering
    \scalebox{0.7}{
    \begin{tikzpicture}[
    vertex/.style={circle, fill=black, inner sep=1.5pt},
    edge/.style={->, >={latex[length=2.5mm, width=1.5mm]}, shorten >=2pt, shorten <=1pt}
    ]

    \begin{scope}
        \node[vertex, label=below:$a$] (a) at (0,0) {};
        \node[vertex, label=below:$b$] (b) at (2.5,0) {};
        \node[vertex, label=above:$c$] (c) at (1.25,2) {};

        \node at (1.25,-1) {$G$};

        \path[edge] (a) edge (c);
        \path[edge] (b) edge (a);
        \path[edge] (c) edge (b);
    \end{scope}

    \begin{scope}[xshift=4.5cm]
        \node[vertex, label=below:$a$] (a) at (0,0) {};
        \node[vertex, label=below:$b$] (b) at (2,0) {};
        \node[vertex, label=above:$c$] (c) at (0,2) {};
        \node[vertex, label=above:$d$] (d) at (2,2) {};

        \node at (1,-1) {$H$};

        \path[edge] (c) edge (a);
        \path[edge] (a) edge (b);
        \path[edge] (d) edge (b);
        \path[edge] (d) edge (c);
        \path[edge] (c) edge (b);
        \path[edge] (a) edge (d);
    \end{scope}

    \begin{scope}[xshift = 9cm]
        \node[vertex, label=below:$a$] (a) at (0,0) {};
        \node[vertex, label=above:$b$] (b) at (0,2) {};
        \node[vertex, label=below:$c$] (c) at (2,0) {};
        \node[vertex, label=above:$d$] (d) at (2,2) {};
        \node[vertex, label=below:$e$] (e) at (4,0) {};
        \node[vertex, label=above:$f$] (f) at (4,2) {};
        \node[vertex, label=below:$g$] (g) at (6,0) {};
        \node[vertex, label=above:$h$] (h) at (6,2) {};

        \node at (3,-1.2) {$T$};

        \path[edge] (a) edge (b);
        \path[edge] (c) edge (a);
        \path[edge] (d) edge (a);
        \path[edge] (e) edge [bend left=20] (a);
        \path[edge] (a) edge (f);
        \path[edge] (g) edge [bend left] (a);
        \path[edge] (h) edge (a);
        \path[edge] (b) edge (c);
        \path[edge] (d) edge (b);
        \path[edge] (e) edge (b);
        \path[edge] (f) edge [bend right=20] (b);
        \path[edge] (b) edge (g);
        \path[edge] (h) edge [bend right] (b);
        \path[edge] (d) edge (c);
        \path[edge] (c) edge (e);
        \path[edge] (f) edge (c);
        \path[edge] (g) edge [bend left=20] (c);
        \path[edge] (h) edge (c);
        \path[edge] (d) edge (f);
        \path[edge] (d) edge (e);
        \path[edge] (d) edge (g);
        \path[edge] (h) edge [bend right=20] (d);
        \path[edge] (f) edge (e);
        \path[edge] (e) edge (g);
        \path[edge] (e) edge (h);
        \path[edge] (f) edge (h);
        \path[edge] (g) edge (f);
        \path[edge] (g) edge (h);
    \end{scope}

    \end{tikzpicture}
    }
    \caption{Tournaments $G$, $H$, and $T$, each with dichromatic number 2, for which a 0-step increases the dichromatic number by 0, 1, and 2, respectively.}
    \label{fig:col ex}
\end{figure}

For a 0-step, we found examples of tournaments for which the dichromatic number remains the same, increases by 1, or increases by 2. See Figure \ref{fig:col ex}. 

The dichromatic number of the base graph, $\vec\chi(G_0)$, gives a trivial lower bound on the dichromatic number of $G_1$, where $G_1$ is obtained via a 0-step. For every $n\geq 3$, there exist oriented graphs on $n$ nodes attaining this lower bound, as the following theorem shows. 

\begin{theorem}
    For every $n\geq 3$, there exist at least $n-2$ non-isomorphic tournaments $G_0$ on $n$ nodes such that $\vec\chi(G_0) = \vec\chi(G_1) = 2$, where $G_1$ is obtained via a 0-step. 
\end{theorem}

\begin{proof}
    For a fixed $n\geq 3$, construct the tournament $G_0$ as follows. Let the nodes be labeled $\{v_1,\dots,v_n\}$. Add the arcs $v_1v_2$, $v_2v_i$ for exactly one $i\in \{3,\dots,n\}$, and $v_jv_k$ for $j>k$, where $(j,k)\not\in \{(2,1),(i,2)\}$. 

    Note that $G_0$ is not transitive, so $\vec\chi(G_0)>1$. The following two partitions each define a 2-dicoloring of $G_0$:
    $$C_1 = \{v_i\}, C_2 = \{v_1,v_2,\dots,v_{i-1},v_{i+1},\dots,v_n\}$$
    $$C_1 = \{v_2\}, C_2 = \{v_1,v_3,v_4,\dots,v_n\}$$

    The following are color classes of a 2-dicoloring of $G_1$:
    $$C_1=\{v_2,v_1',v_2',\dots,v_{i-1}',v_{i+1}',\dots,v_n'\}, C_2=\{v_i',v_1,v_3,v_4,\dots,v_n\}.$$
    Since the dichromatic number of $G_1$ is at least the dichromatic number of $G_0$, this shows that $\vec\chi(G_1) = \vec\chi(G_0) = 2$. For each choice of $i$, the unique strongly connected component of $G_0$ with more than one node is $\{v_1,\dots,v_i\}$. Its order distinguishes the constructions for different values of $i$, giving $n-2$ pairwise non-isomorphic tournaments.
\end{proof}

In general, the dichromatic number can at most double as a result of a 0-step. For $t\geq 1$ and $2t$ consecutive 0-steps, we obtain an upper bound on the dichromatic number of $G_{2t}$ using a partition of $\mathbb{Z}_2^{2t}$. 

For $k\geq 1$, define the function $f_k:\mathbb{Z}_2^k \times V(G_0) \to V(G_k)$ recursively as follows. For $k=1$, let $f_1(0,v)=v$ and $f_1(1,v) = v'$, the clone of $v$ in $G_1$, for all $v\in V(G_0)$. For $k>1$, we have 
$$f_k(b_1b_2\dots b_{k-1}b_k,v) = \begin{cases}
    f_{k-1}(b_1b_2\dots b_{k-1},v),& \text{if~} b_k=0\\ 
    f_{k-1}(b_1b_2\dots b_{k-1},v)', & \text{if~}b_k=1
\end{cases}$$

Now, for $t\geq 1$, let $h:\mathbb{Z}_2^{2t}\times V(G_0) \to V(G_{2t})$ be the function $h(B,v) = f_{2t}(B,v)$ for every node $v\in V(G_0)$ and binary string $B\in \mathbb{Z}_2^{2t}$. It is straightforward to show that $h$ is a bijection. We have the following lemma that describes how the strings $B_1$ and $B_2$ corresponding to nodes $h(B_1,u)$ and $h(B_2,v)$ in $G_{2t}$ determine the direction of the arc between the two nodes, if such an arc exists. 

\begin{lemma}
    Let $G_0$ be an oriented graph, $t\geq 1$, and $G_{2t}$ the result of $2t$ consecutive 0-steps. Let $h:\mathbb{Z}_2^{2t}\times V(G_0) \to V(G_{2t})$ be defined as above, and $B_1,B_2 \in \mathbb{Z}_2^{2t}$, not necessarily distinct. If $(u,v)$ is an arc in $G_0$ and $B_1 \cdot B_2 \equiv 0 \pmod{2}$, then $(h(B_1,u),h(B_2,v))$ is an arc in $G_{2t}$. If, instead, $B_1 \cdot B_2 \equiv 1 \pmod{2}$, then $(h(B_2,v),h(B_1,u))$ is an arc in $G_{2t}$. 
    \label{lem: strings}
\end{lemma}

\begin{proof}
    We use induction on $t$. First, let $t=1$. If $B_1,B_2\in \mathbb{Z}^2_2$, then we can consider each of the cases for $B_1$ and $B_2$ and the corresponding arc in $G_2$. In each case, we consider whether $h(B_1,u)$ is a clone or parent node in $G_2$, and whether it is the clone or parent of $u$ or $u'$, the clone of $u$ in $G_1$, and similarly for $h(B_2,v)$. The arguments for each case are summarized in Table \ref{tab:t1}, with the case $B_1 \cdot B_2 \equiv 0 \pmod{2}$ in the first section of the table, and $B_1 \cdot B_2 \equiv 1 \pmod{2}$ in the second section. 

    \begin{table}[htbp]
    \centering
    \begin{tabular}{c@{~~~}c@{~}|@{~~~}c@{~~~}|@{~~~}c@{~~~}|@{~~~}c@{~~~}}
        $B_1$ & $B_2$ & $G_0$ & $G_1$ & $G_2$ \\
        \hline
        00 & 00 & $(u,v)$ & $(u,v)$ & $(u,v) = (h(00,u),h(00,v))$\\
        00 & 01 & $(u,v)$ & $(u,v)$ & $(u,{\bf v''})= (h(00,u),h(01,v))$\\
        00 & 10 & $(u,v)$ & $(u,{\bf v'})$ & $(u,v')= (h(00,u),h(10,v))$\\
        00 & 11 & $(u,v)$ & $(u,{\bf v'})$ & $(u,{\bf (v')'})= (h(00,u),h(11,v))$\\
        01 & 00 & $(u,v)$ & $(u,v)$ & $({\bf u''},v)= (h(01,u),h(00,v))$\\
        01 & 10 & $(u,v)$ & $(u,{\bf v'})$ & $({\bf u''},v')= (h(01,u),h(10,v))$\\
        10 & 00 & $(u,v)$ & $({\bf u'},v)$ & $(u',v)= (h(10,u),h(00,v))$\\
        10 & 01 & $(u,v)$ & $({\bf u'},v)$ & $(u',{\bf v''})= (h(10,u),h(01,v))$\\
        11 & 00 & $(u,v)$ & $({\bf u'},v)$ & $({\bf (u')'},v)= (h(11,u),h(00,v))$\\
        11 & 11 & $(u,v)$ & $({\bf v'},{\bf u'})$ & $({\bf (u')'},{\bf(v')'})= (h(11,u),h(11,v))$\\
        \hline
        01 & 01 & $(u,v)$ & $(u,v)$ & $({\bf v''},{\bf u''}) = (h(01,v),h(01,u))$\\
        01 & 11 & $(u,v)$ & $(u,{\bf v'})$ & $({\bf (v')'},{\bf u''})= (h(11,v),h(01,u))$\\
        10 & 10 & $(u,v)$ & $({\bf v'},{\bf u'})$ & $(v',u')= (h(10,v),h(10,u))$\\
        10 & 11 & $(u,v)$ & $({\bf v'},{\bf u'})$ & $({\bf (v')'},u')= (h(11,v),h(10,u))$\\
        11 & 01 & $(u,v)$ & $({\bf u'},v)$ & $({\bf v''},{\bf (u')'})= (h(01,v),h(11,u))$\\
        11 & 10 & $(u,v)$ & $({\bf v'},{\bf u'})$ & $(v',{\bf (u')'})= (h(10,v),h(11,u))$\\
    \end{tabular}
    \caption{A summary of the $t=1$ case, with the corresponding arcs in $G_1$ and $G_2$ listed. We use $u''$ to denote the clone of $u$ in $G_2$ and $(u')'$ to denote the clone of $u'$ in $G_2$. To improve readability, new clones at each time-step are in bold.}
    \label{tab:t1}
\end{table}

    Thus, the result holds for $t=1$. Suppose next that $t>1$ and that the statement holds for $1\leq k<t$. Denote by $h' : \mathbb{Z}^{2(t-1)}_2 \times V(G_0) \to V(G_{2(t-1)})$ the function $h'(B,u) = f_{2(t-1)}(B,u)$. 
    
    First, suppose $B_1, B_2 \in \mathbb{Z}^{2t}_2$ satisfy $B_1 \cdot B_2 \equiv 0 \pmod{2}$ and consider an arc $(u,v)$ in $G_0$. 
    The string $B_1$ has the form $B_1 = B_1'a_1a_2$, where $B_1'\in \mathbb{Z}^{2(t-1)}_2$. Similarly, $B_2$ has the form $B_2 = B_2'c_1c_2$, where $B_2'\in \mathbb{Z}^{2(t-1)}_2$. We consider the direction of the arc between nodes $h(B_1,u)$ and $h(B_2,v)$ for each case of the values of $a_1,a_2,c_1,c_2$. 

    Suppose that $a_1a_2 = 00$ and $c_1c_2 = 00$. It must be the case that $B_1' \cdot B_2' \equiv 0 \pmod{2}$. By induction, we have the arc $(h'(B'_1,u),h'(B'_2,v))$ in $G_{2(t-1)}$. The same is true in each subsequent case where $a_1a_2 \cdot c_1c_2 \equiv 0 \pmod{2}$. Now, by definition of $h$, $h(B'_100,u) =h'(B'_1,u)$ and, similarly, $h(B'_200,v) = h'(B'_2,v)$. Therefore, we have the arc $(h(B_1,u),h(B_2,v))$ in $G_{2(t-1)}$ and thus also in $G_{2t}$.

    Next, suppose that $a_1a_2 = 01$ and $c_1c_2 = 01$. It must be the case that $B_1' \cdot B_2' \equiv 1 \pmod{2}$. By induction, we have the arc $(h'(B'_2,v),h'(B'_1,u))$ in $G_{2(t-1)}$. The same is true in each subsequent case where $a_1a_2 \cdot c_1c_2 \equiv 1 \pmod{2}$. By definition of $h$, the node $h(B_1,u) = h(B'_101,u)$ is the clone of the node $h'(B'_1,u)$ in $G_{2t}$. Similarly, the node $h(B_2,v)$ is the clone of $h'(B'_2,v)$ in $G_{2t}$. Since $(h'(B'_2,v),h'(B'_1,u))$ is an arc in $G_{2(t-1)}$, it is also an arc in $G_{2t-1}$, and, cloning both nodes, $(h(B_1,u),h(B_2,v))$ is an arc in $G_{2t}$.

    The remaining cases are summarized in Table~\ref{tab:t2}. 

    \begin{table}[htbp]
        \centering
        \begin{tabular}{c@{~~~}c@{~}|@{~~~}c@{~~~}|@{~~~}c@{~~~}|@{~~~}c@{~~~}}
           $a_1a_2$ & $c_1c_2$ & $G_{2(t-1)}$ & $G_{2t-1}$ & $G_{2t}$  \\
           \hline
            00 & 00 & $(h'(B_1',u),h'(B_2',v))$ & $(h'(B_1',u),h'(B_2',v))$ & $(h(B_1'00,u),h(B_2'00,v))$\\
            00 & 01 & $(h'(B_1',u),h'(B_2',v))$ & $(h'(B_1',u),h'(B_2',v))$ & $(h(B_1'00,u),{\bf h(B_2'01,v)})$\\
            00 & 10 & $(h'(B_1',u),h'(B_2',v))$ & $(h'(B_1',u),{\bf h'(B_2',v)'})$ & $(h(B_1'00,u),h(B_2'10,v))$\\
            00 & 11 & $(h'(B_1',u),h'(B_2',v))$ & $(h'(B_1',u),{\bf h'(B_2',v)'})$ & $(h(B_1'00,u),{\bf h(B_2'11,v)})$\\
            01 & 00 & $(h'(B_1',u),h'(B_2',v))$ & $(h'(B_1',u),h'(B_2',v))$ & $({\bf h(B_1'01,u)},h(B_2'00,v))$\\
            01 & 10 & $(h'(B_1',u),h'(B_2',v))$ & $(h'(B_1',u),{\bf h'(B_2',v)'})$ & $({\bf h(B_1'01,u)},h(B_2'10,v))$\\
            10 & 00 & $(h'(B_1',u),h'(B_2',v))$ & $({\bf h'(B_1',u)'},h'(B_2',v))$ & $(h(B_1'10,u),h(B_2'00,v))$\\
            10 & 01 & $(h'(B_1',u),h'(B_2',v))$ & $({\bf h'(B_1',u)'},h'(B_2',v))$ & $(h(B_1'10,u),{\bf h(B_2'01,v)})$\\
            11 & 00 & $(h'(B_1',u),h'(B_2',v))$ & $({\bf h'(B_1',u)'},h'(B_2',v))$ & $({\bf h(B_1'11,u)},h(B_2'00,v))$\\
            11 & 11 & $(h'(B_1',u),h'(B_2',v))$ & $({\bf h'(B_2',v)'},{\bf h'(B_1',u)'})$ & $({\bf h(B_1'11,u)},{\bf h(B_2'11,v)})$\\
            \hline
            01 & 01 & $(h'(B_2',v),h'(B_1',u))$ & $(h'(B_2',v),h'(B_1',u))$ & $({\bf h(B_1'01,u)},{\bf h(B_2'01,v)})$\\
            01 & 11 & $(h'(B_2',v),h'(B_1',u))$ & $({\bf h'(B_2',v)'},h'(B_1',u))$ & $({\bf h(B_1'01,u)},{\bf h(B_2'11,v)})$\\
            10 & 10 & $(h'(B_2',v),h'(B_1',u))$ & $({\bf h'(B_1',u)'},{\bf h'(B_2',v)'})$ & $(h(B_1'10,u),h(B_2'10,v))$\\
            10 & 11 & $(h'(B_2',v),h'(B_1',u))$ & $({\bf h'(B_1',u)'},{\bf h'(B_2',v)'})$ & $(h(B_1'10,u),{\bf h(B_2'11,v)})$\\
            11 & 01 & $(h'(B_2',v),h'(B_1',u))$ & $({\bf h'(B_2',v)'},h'(B_1',u))$ & $({\bf h(B_1'11,u)},{\bf h(B_2'01,v)})$\\
            11 & 10 & $(h'(B_2',v),h'(B_1',u))$ & $({\bf h'(B_1',u)'},{\bf h'(B_2',v)'})$ & $({\bf h(B_1'11,u)},h(B_2'10,v))$\\
        \end{tabular}
        \caption{A summary of the $t>1$ case. Clone nodes at each time-step are bolded.}
        \label{tab:t2}
    \end{table}

Therefore, if $B_1 \cdot B_2 \equiv 0 \pmod{2}$ and $(u,v)$ is an arc in $G_0$, then $(h(B_1,u),h(B_2,v))$ is an arc in $G_{2t}$. 

Next, suppose $B_1, B_2 \in \mathbb{Z}^{2t}_2$ satisfy $B_1 \cdot B_2 \equiv 1 \pmod{2}$ and consider an arc $(u,v)$ in $G_0$. 
    As before, the string $B_1$ has the form $B_1 = B_1'a_1a_2$, where $B_1'\in \mathbb{Z}^{2(t-1)}_2$. Similarly, $B_2$ has the form $B_2 = B_2'c_1c_2$, where $B_2'\in \mathbb{Z}^{2(t-1)}_2$. 

The only change from the above argument is the induction hypothesis. Now, if $a_1a_2 \cdot c_1c_2 \equiv 0 \pmod{2}$ then it must be that $B'_1\cdot B'_2 \equiv 1 \pmod{2}$, and so $(h'(B'_2,v),h'(B'_1,u))$ is an arc in $G_{2(t-1)}$, by the induction hypothesis.

    Similarly if $a_1a_2 \cdot c_1c_2 \equiv 1 \pmod{2}$, then it must be that $B'_1 \cdot B'_2 \equiv 0 \pmod{2}$, and so $(h'(B'_1,u),h'(B'_2,v))$ is an arc in $G_{2(t-1)}$, by the induction hypothesis. Thus, we can swap the roles of $h'(B'_1,u)$ and $h'(B'_2,v)$ in each of the above arguments and obtain the desired result. 
\end{proof}

Now, given a dicoloring of the oriented graph $G_0$, Lemma~\ref{lem: strings} gives a tool for creating a dicoloring of $G_{2t}$ using the color classes from $G_0$, where $G_{2t}$ is obtained from $G_0$ by $2t$ consecutive 0-steps. 

\begin{theorem} \label{thm:upper coloring}
    Let $G_0$ be an oriented graph and, for $t\geq 0$, let $G_{2t}$ be the result of $2t$ consecutive 0-steps. We then have that
    $$\vec\chi(G_{2t}) \leq (2^{t+1}-1)\vec\chi(G_0).$$
\end{theorem}

\begin{proof} The case $t=0$ is immediate, so assume that $t\geq1$. Define the function $h:\mathbb{Z}^{2t}_2\times V(G_0) \to V(G_{2t})$ as above and, for $S\subseteq V(G_0)$ and $\mathcal{B}\subseteq \mathbb{Z}^{2t}_2$ denote by $(\mathcal{B},S)$ the set $\{h(B,v) \mid B\in \mathcal{B}, v\in S\}$. All congruences are modulo 2. 

    Let $k=2^{t+1}-1$. We will first recursively construct a partition of $\mathbb{Z}^{2t}_2$, $\mathscr{B}_t = \{\mathcal{B}_1, \mathcal{B}_2,\dots,\mathcal{B}_k\}$ and then argue that if $S\subseteq V(G_0)$ induces an acyclic oriented subgraph of $G_0$, then $(\mathcal{B}_i,S)$ induces an acyclic oriented subgraph of $G_{2t}$, for $i=1,2,\dots,k$. 

    Each set $\mathcal{B}_i$ in $\mathscr{B}_t$ will satisfy one of the following two properties:
    \begin{equation} \label{eq:A1}
        B_1 \cdot B_2 \equiv 0 \quad \text{for all~} B_1,B_2\in \mathcal{B}_i \tag{A1}
    \end{equation}
    \begin{equation} \label{eq:A2}
        B_1 \cdot B_2 \equiv 1 \quad \text{for all~} B_1,B_2\in \mathcal{B}_i \tag{A2}
    \end{equation}
    where in each case, the binary strings $B_1$ and $B_2$ are not necessarily distinct. 

    For $\mathcal{B}\subseteq \mathbb{Z}^{2(t-1)}_2$, let $\mathcal{B} \oplus a_1a_2$ denote the set $\{b_1b_2\dots b_{2(t-1)}a_1a_2 \mid b_1b_2\dots b_{2(t-1)}\in \mathcal{B}\}$. We construct $\mathscr{B}_t$ recursively as follows. For $t=1$, let $\mathscr{B}_1$ be the partition $\mathcal{B}_1 = \{00,11\}$, $\mathcal{B}_2 = \{10\}$, $\mathcal{B}_3 = \{01\}$. 

    We will show that the following property is preserved by the construction. For $t\geq 1$, let $\mathcal{B}_1 = \{b_1b_2\dots b_{2t} \mid b_{2i-1}=b_{2i}, i=1,2,\dots,t\}$. \begin{itemize}
        \item[($*$)] For each set $\mathcal{B} \neq \mathcal{B}_1$ in $\mathscr{B}_t$ that satisfies (\ref{eq:A1}), there is a set $\mathcal{B}^* \in \mathscr{B}_t$ that \begin{itemize}[leftmargin=.5in]
            \item[(i)] also satisfies (\ref{eq:A1}), and 
            \item [(ii)] satisfies $B_1 \cdot B_2 \equiv 1$ for every $B_1 \in \mathcal{B}$ and $B_2\in \mathcal{B}^*$.
            \end{itemize} 
            Similarly, for each set $\mathcal{B}$ in the partition that satisfies (\ref{eq:A2}), there is a set $\mathcal{B}^*$ in the partition that \begin{itemize}[leftmargin=.5in]
                \item[(i)] also satisfies (\ref{eq:A2}), and 
                \item[(ii)] satisfies $B_1 \cdot B_2 \equiv 0$ for every $B_1 \in \mathcal{B}$ and $B_2\in \mathcal{B}^*$.
                \end{itemize}
    \end{itemize}

We also require that the parts other than $\mathcal{B}_1$ are partitioned into pairs satisfying $(*)$. Thus, each such part
$\mathcal{B}$ has a designated partner $\mathcal{B}^*\neq\mathcal{B}$, with $(\mathcal{B}^*)^*=\mathcal{B}$.
For $t=1$, pair $\{10\}$ with $\{01\}$; these parts satisfy $(*)$ since both satisfy (\ref{eq:A2}) and $10\cdot01\equiv0$.
At each subsequent stage, we use these designated partners in the construction below. Now, for $t>1$, let $\mathscr{B}_{t-1} = \{\mathcal{B}'_1,\mathcal{B}'_2,\dots,\mathcal{B}'_{\frac{k-1}{2}}\}$.  

    For $i=1,2,\dots,\frac{k-1}{2}$, let $\mathcal{B}_i = \mathcal{B}'_i \oplus 00 \cup \mathcal{B}'_i \oplus 11$. If $\mathcal{B}'_i$ satisfies condition (\ref{eq:A1}), then $\mathcal{B}_i$ does as well. Further, if $i\neq1$ and $\mathcal{B}'_i$ satisfies (\ref{eq:A1}), then there is a set $\mathcal{B}'_j=\mathcal{B}'^*_i$ satisfying property ($*$). The sets $\mathcal{B}_i$ and $\mathcal{B}_j = \mathcal{B}^*_i$, constructed from this pair, also satisfy ($*$).  Similarly for condition (\ref{eq:A2}). 

    Let $\mathcal{B}_{\frac{k+1}{2}} = \mathcal{B}'_1 \oplus 01$ and $\mathcal{B}_{\frac{k+3}{2}} = \mathcal{B}'_1 \oplus 10$. Since $\mathcal{B}'_1$ satisfies (\ref{eq:A1}), both $\mathcal{B}_{\frac{k+1}{2}}$ and $\mathcal{B}_{\frac{k+3}{2}}$ satisfy (\ref{eq:A2}). Note that $B_1 \cdot B_2 \equiv 0$ for $B_1 \in \mathcal{B}_{\frac{k+1}{2}}$ and $B_2 \in \mathcal{B}_{\frac{k+3}{2}}$, satisfying ($*$). 

    Now, for $i=2,\dots,\frac{k-1}{2}$, let $\mathcal{B}_{i+\frac{k+1}{2}} = \mathcal{B}'_i \oplus 10 \cup \mathcal{B}'^*_i \oplus 01$. If both $\mathcal{B}'_i$ and $\mathcal{B}'^*_i$ satisfy (\ref{eq:A1}), then $\mathcal{B}_{i+\frac{k+1}{2}}$ satisfies (\ref{eq:A2}). Similarly, if both satisfy (\ref{eq:A2}), then $\mathcal{B}_{i+\frac{k+1}{2}}$ satisfies (\ref{eq:A1}). Moreover, if $\mathcal{B}'^*_i = \mathcal{B}'_j$, then $\mathcal{B}_{i+\frac{k+1}{2}} = \mathcal{B}'_i \oplus 10 \cup \mathcal{B}'_j \oplus 01$ and $\mathcal{B}_{j+\frac{k+1}{2}} = \mathcal{B}'_i\oplus 01 \cup \mathcal{B}'_j \oplus 10$. This pair of sets satisfies property ($*$). 
    
    Note that there are a total of $\frac{k+1}{2}+\frac{k-1}{2} = k$ sets $\mathcal{B} \subseteq \mathbb{Z}^{2t}_2$ formed by this construction. It is straightforward to show that these sets form a partition $\mathscr{B}_t$ of $\mathbb{Z}^{2t}_2$. Further, each set satisfies either (\ref{eq:A1}) or (\ref{eq:A2}), and $\mathscr{B}_t$ satisfies ($*$).

    We now show that the partition $\mathscr{B}_t$ can be used to construct a $k\vec\chi(G_0)$-dicoloring of $G_{2t}$. 
    Let $S_1,S_2,\dots,S_{\vec\chi(G_0)} \subseteq V(G_0)$ be the color classes of a $\vec\chi(G_0)$-dicoloring of $G_0$. We show that, for $i=1,2,\dots,k$ and $j=1,2,\dots,\vec\chi(G_0)$, the oriented subgraph of $G_{2t}$ induced by $(\mathcal{B}_i,S_j)$ is acyclic, by considering the two cases for $\mathcal{B}_i$. \\

\noindent \emph{Case 1:} Suppose that $\mathcal{B}_i$ satisfies condition (\ref{eq:A1}). Label the nodes of $S_j$ as $v_1,v_2,\dots,v_r$ so that every arc $(v_a,v_b)$ with both endpoints in $S_j$ satisfies $a<b$. We will show that the binary strings in $\mathcal{B}_i$ can be ordered as $B_1,\dots,B_m$ so that, in the node ordering
$$h(B_1,v_1),\dots,h(B_m,v_1),h(B_1,v_2),\dots,h(B_m,v_2),\dots,h(B_1,v_r),\dots,h(B_m,v_r),$$
every arc between listed nodes points from an earlier node to a later node. Consequently, these nodes induce an acyclic oriented subgraph of $G_{2t}$.

By Lemma~\ref{lem: strings}, if $(v_a,v_b)$ is an arc in $G_0$, then $(h(B,v_a),h(B',v_b))$ is an arc in $G_{2t}$ for any $B,B'\in\mathcal{B}_i$. It remains to be shown that $h(B_1,v_a),\dots,h(B_m,v_a)$ induces an acyclic oriented graph, for $a=1,2,\dots,r$.

Each node $h(B,v_a)$ with $B\in\mathcal{B}_i$ projects onto $v_a$ in $G_0$, using the projection function $p:V(G_{2t})\to V(G_0)$ defined in Section~2. By Lemma~\ref{lem:projection}, the nodes $h(B_1,v_a),\dots,h(B_m,v_a)$ induce a tournament in $G_{2t}$. Thus, it suffices to show that the subgraph induced by this set of nodes contains no directed 3-cycle.

For two distinct strings with longest common prefix $A$, write them as $A\,1\,C_1$ and $A\,0\,C_2$. At the first differing coordinate, the corresponding nodes are a clone and its parent, so the arc points from the former to the latter. At subsequent steps, this arc reverses exactly when both nodes are cloned. Hence, the final arc points from $h(A\,1\,C_1,v_a)$ to $h(A\,0\,C_2,v_a)$ if $C_1\cdot C_2\equiv0$, and in the opposite direction if $C_1\cdot C_2\equiv1$. Here, the dot product of two empty suffixes is zero. Suppose we have arcs $(h(B_1,v_a),h(B_2,v_a))$ and $(h(B_2,v_a),h(B_3,v_a))$. The strings $B_1$ and $B_2$ have the form
$$B_1 = A\,1\,C_1 \quad \text{and} \quad B_2 = A\,0\,C_2$$
with $C_1\cdot C_2 \equiv 0$, or $C_1,C_2 = \emptyset$, or
$$B_1 = A\,0\,C_1 \quad \text{and} \quad B_2 = A\,1\,C_2$$
with $C_1 \cdot C_2 \equiv 1$.

Similarly, the strings $B_2$ and $B_3$ have the form
$$B_2 = D\,1\,C_3 \quad \text{and} \quad B_3 = D\,0\,C_4$$
with $C_3\cdot C_4 \equiv 0$, or $C_3,C_4 = \emptyset$, or
$$B_2 = D\,0\,C_3 \quad \text{and} \quad B_3 = D\,1\,C_4$$
with $C_3 \cdot C_4 \equiv 1$.

Using this information, we consider the arc between nodes $h(B_1,v_a)$ and $h(B_3,v_a)$. First, suppose that $|D|<|A|$. If $C_3\cdot C_4 \equiv 0$, then $B_2 = D\,1\,C_3$ and we can write $A = D\,1\,F$. Since $B_2 \cdot B_3 \equiv 0$, we must have $D\cdot D \equiv 0$.

If $B_1 = A\,1\,C_1$, then we can rewrite it as $B_1 = D\,1\,F\,1\,C_1$. Recall that $B_3 = D\,0\,C_4$. We know that $B_1 \cdot B_3 \equiv 0$ and $D\cdot D \equiv 0$ so it must be the case that $(F\,1\,C_1)\cdot C_4 \equiv 0$ as well.

The case $B_1 = A\,0\,C_1$ is similar. Thus, if $C_3 \cdot C_4 \equiv 0$, then there is an arc from $h(B_1,v_a)$ to $h(B_3,v_a)$ in $G_{2t}$.

If $C_3 \cdot C_4 \equiv 1$, then $B_2 = D\,0\,C_3$. As above, we write $A = D\,0\,F$ and, from $B_2 \cdot B_3 \equiv 0$, conclude that $D \cdot D \equiv 1$. If $B_1 = A\,1\,C_1$, then we can rewrite it as $B_1 = D\,0\,F\,1\,C_1$. Recall that $B_3 = D\,1\,C_4$. Since $B_1 \cdot B_3 \equiv 0$ and $D \cdot D \equiv 1$, we have that $(F\,1\,C_1)\cdot C_4 \equiv 1$.

The case $B_1 = A\,0\,C_1$ is similar. Thus, if $C_3 \cdot C_4 \equiv 1$, then there is an arc from $h(B_1,v_a)$ to $h(B_3,v_a)$ in $G_{2t}$. A similar argument gives the same arc when $|D|>|A|$.

Finally, suppose that $|D|=|A|$. If $C_3 \cdot C_4 \equiv 0$, then $D \cdot D \equiv 0$ and $B_1$ has the form $B_1 = D\,0\,C_1$. Since $(h(B_1,v_a),h(B_2,v_a))$ is an arc, we must have $C_1\cdot C_3 \equiv 1$ and so $D \cdot D \equiv 1$, a contradiction. Similarly, the case of $C_3 \cdot C_4 \equiv 1$ results in a contradiction.

Therefore, $\{h(B_1,v_a),h(B_2,v_a),\dots,h(B_m,v_a)\}$ induces an acyclic tournament in $G_{2t}$, for $a=1,2,\dots,r$. This shows that $(\mathcal{B}_i,S_j)$ induces an acyclic oriented subgraph of $G_{2t}$ for $j=1,2,\dots,\vec\chi(G_0)$ and $\mathcal{B}_i$ satisfying (\ref{eq:A1}). \\

\noindent \emph{Case 2:} Suppose that $\mathcal{B}_i$ satisfies condition (\ref{eq:A2}). Label the nodes of $S_j$ as $v_1,v_2,\dots,v_r$ so that every arc $(v_a,v_b)$ with both endpoints in $S_j$ satisfies $a<b$. By Lemma~\ref{lem: strings} and an argument analogous to Case 1, the binary strings in $\mathcal{B}_i$ can be ordered as $B_1,\dots,B_m$ so that, in the node ordering
$$h(B_1,v_r),\dots,h(B_m,v_r),h(B_1,v_{r-1}),\dots,h(B_m,v_{r-1}),\dots,h(B_1,v_1),\dots,h(B_m,v_1),$$
every arc between listed nodes points from an earlier node to a later node. Thus, $(\mathcal{B}_i,S_j)$ induces an acyclic oriented subgraph of $G_{2t}$ for every $j=1,2,\dots,\vec\chi(G_0)$ and every $\mathcal{B}_i$ satisfying (\ref{eq:A2}).

Since the sets $(\mathcal{B}_i,S_j)$ partition the nodes of $G_{2t}$ and each induces an acyclic oriented subgraph of $G_{2t}$, they form the color classes of a $k\vec\chi(G_0)$-dicoloring of $G_{2t}$, where $k=2^{t+1}-1$.
\end{proof}

To obtain a lower bound on the dichromatic number, we consider the tournament formed from a single node by successive 0-steps. Double counting directed paths of length two bounds the cardinality of its acyclic sets, and hence, the order of each color class.

\begin{theorem}\label{thm:lower coloring}
Let $G_0$ be an oriented graph, and let $G_t$ be obtained from $G_0$ via $t$ consecutive 0-steps. We then have that
$$\max\left\{\vec\chi(G_0),\left\lceil\frac{2^t}{\sqrt{3\cdot2^t-2}}\right\rceil\right\}\leq \vec\chi(G_t).$$
\end{theorem}

\begin{proof}
Since $G_0$ is an induced subgraph of $G_t$, we have that $\vec\chi(G_0)\leq\vec\chi(G_t)$. The remaining bound is immediate if $t=0$, so assume $t\geq1$. Let $F_0$ consist of a single node of $G_0$. For $j=0,1,\ldots,t-1$, let $F_{j+1}$ be the subgraph induced by the nodes of $F_j$ and their clones at the next time-step. Thus, $F_j$ is a tournament on $2^j$ nodes and is an induced subgraph of $G_j$. 
Equivalently, if $v$ is the node of $F_0$, then for $j\geq1$, the graph $F_j$ is induced by the set $\{x\in V(G_j):p(x)=v\}$, where $p$ denotes the projection onto $G_0$. Thus, the recursive construction collects all nodes whose projection is $v$.

We first count directed paths of length two with prescribed endpoints. Throughout the remainder of the proof, all paths mentioned are directed and have length two.

We claim that any two distinct nodes of $F_j$ are the endpoints, in either order, of exactly $(2^j-2)/2=2^{j-1}-1$ such paths. This will allow us to count paths whose endpoints belong to an acyclic set.

We proceed by induction on $j\geq1$. The claim holds for $F_1$, which has two nodes and no directed path of length two. Suppose it holds for $F_j$, and consider two distinct nodes $x,y$ of $F_{j+1}$. Since we count paths with these endpoints in either order, we may interchange $x$ and $y$. We consider three cases.

If $x=u$ is a parent and $y=u'$ is its clone, then every other clone is the middle node of one such path, and no parent is. This gives $2^j-1$ paths.

If $x,y$ are both parents $u,v$, or both clones $u',v'$, then label the endpoints so that $(u,v)$ is an arc. Each middle node counted for $u,v$ in $F_j$ and its clone are middle nodes of paths with endpoints $x,y$ in $F_{j+1}$, giving two paths for each middle node in $F_j$. There is one additional middle node: $v'$ when $x=u$ and $y=v$, and $u$ when $x=u'$ and $y=v'$. The other unused node among $u,u',v,v'$ gives no additional path. Thus, the count is $2(2^j-2)/2+1=2^j-1$.

Finally, suppose $x=u$ and $y=v'$, where $u\neq v$. For each $z\in V(F_j)\setminus\{u,v\}$, exactly one of $z,z'$ is a middle node, since cloning $z$ preserves its arc with $x$ and reverses its arc with $y$. Exactly one of $u',v$ also contributes: $u'$ if $(u,v)$ is an arc, and $v$ otherwise. In this case, there are also $2^j-1$ paths. Thus, every pair of distinct endpoints in $F_{j+1}$ has exactly $2^j-1=(2^{j+1}-2)/2$ paths, proving the claim.

We now use this count to bound the cardinality of an acyclic set $S$ in $F_t$. Let $m=|S|$ and suppose $m>0$. We count paths with both endpoints in $S$. Counting by their unordered pairs of endpoints gives exactly $\binom{m}{2}(2^t-2)/2$ paths.

We can also count these paths by their middle nodes. If the middle node lies in $S$, then all three nodes belong to $S$. Since $S$ induces a transitive tournament, each three-node subset contains exactly one directed path of length two. Thus, middle nodes in $S$ contribute exactly $\binom{m}{3}$ paths.

Now fix a node $z\notin S$, and suppose it has $a$ in-neighbors in $S$. Since $F_t$ is a tournament, each of the $m$ nodes in $S$ is either an in-neighbor or an out-neighbor of $z$, so $z$ has $m-a$ out-neighbors in $S$. Choosing one in-neighbor as the first node and one out-neighbor as the last node gives exactly $a(m-a)$ paths with middle node $z$. This product is largest when the $m$ nodes are divided as evenly as possible between the two groups, giving at most $(m/2)(m/2)=m^2/4$ paths. There are $2^t-m$ nodes outside $S$, so these middle nodes contribute at most $(2^t-m)m^2/4$ paths.

Hence, there are at most $\binom{m}{3}+(2^t-m)\frac{m^2}{4}$ paths with both endpoints in $S$. Comparing this upper bound with the exact count by endpoints gives
$$\binom{m}{2}\frac{2^t-2}{2}\leq\binom{m}{3}+(2^t-m)\frac{m^2}{4}.$$
We therefore have that $m(m^2-3\cdot2^t+2)\leq0$. Since $m>0$, we obtain $m^2\leq3\cdot2^t-2$.

Every color class in a dicoloring of $F_t$ therefore contains at most $\sqrt{3\cdot2^t-2}$ nodes. Since $F_t$ has $2^t$ nodes, it requires at least $\lceil2^t/\sqrt{3\cdot2^t-2}\rceil$ colors. As $F_t$ is an induced subgraph of $G_t$, the result follows.
\end{proof}

Combining Theorems~\ref{thm:upper coloring} and ~\ref{thm:lower coloring}, we obtain, after $2t$ consecutive 0-steps,
$$\frac{2^t}{\sqrt3}\leq\vec\chi(G_{2t})\leq(2^{t+1}-1)\vec\chi(G_0).$$
Thus, for each fixed initial oriented graph $G_0$, we have that $\vec\chi(G_{2t})=\Theta(2^t)$. 

\section{Conclusion and Future Directions}

The ILTO model generates dense oriented graphs with strong connectivity properties under suitable assumptions on the initial graph and generating sequence. For a generating sequence with infinite support, an oriented graph $H$ eventually appears as an induced subgraph if and only if $N_H^*$ is an induced subgraph of $G_0^*$. When the initial graph is a tournament and the generating sequence has infinite support, the resulting ILTO graphs form a quasirandom sequence of tournaments. We also studied how domination numbers, the cop number, and the dichromatic number behave under a 0- or 1-step.

For a fixed base graph, Theorem~\ref{thm:bounded indomination} shows that the in-domination number is bounded independently of the generating sequence and the number of time-steps. Universality implies that the dichromatic number tends to infinity whenever the generating sequence has infinite support, and our bounds give $\vec\chi(G_{2t})=\Theta(2^t)$ after $2t$ consecutive 0-steps. Determining its asymptotic growth for general generating sequences with infinite support remains open. For each fixed base graph $G_0$, does $\vec\chi(G_{2t})/2^t$ converge as $t\to\infty$ when all steps are 0-steps? If so, what is its limit? Since $\alpha(G_t)=\alpha(G_0)$, Theorem~\ref{thm: cops} and the invariance of the cop number under 1-steps give $c(G_t)\leq3\alpha(G_0)$ after the first 0-step. Determining the exact behavior of the cop number within this bound remains open. 

A probabilistic version of ILTO would apply a 0-step with probability $p$ and a 1-step with probability $1-p$, independently at each time-step. Another variation would omit arcs between clone nodes while retaining the other construction rules. The new clones would then form an independent set, as in the ILM \cite{ILM} and ILDT \cite{ILDT} models. The effects of these variations on densification, connectivity, and universality remain to be investigated.

\end{document}